\documentclass[review]{siamart0516}

\newcommand{\cahiernumber}{27}  % Insert your Cahier du GERAD number.

\makeatletter
\def\cahierline{\gdef\@cahierline}
\cahierline{Cahier du GERAD G-\number\year-\cahiernumber}
\def\ps@firstpage{%\ps@plain
  \def\@oddhead{\hfil\normalfont\small\@cahierline}%
  \let\@evenhead\@oddhead % in case an article starts on a left-hand page
}
\def\ps@myheadings{%
  \def\@oddfoot{\hfil\normalfont\small\color{gray}{\@cahierline}}%
  \let\@evenfoot\@oddfoot % in case an article starts on a left-hand page
  \def\@evenhead{\rlap{\thepage}\hfil\upshape\footnotesize\leftmark\hfil}%          %!
  \def\@oddhead{\hbox{}\hfil{\upshape\footnotesize\rightmark}\hfil\llap{\thepage}}% %!
  \let\@mkboth\@gobbletwo
  \let\sectionmark\@gobble
  \let\subsectionmark\@gobble
}
\makeatother

\usepackage{multirow}
\usepackage[low-sup]{subdepth}
\usepackage{microtype}
\usepackage{graphicx}%, subfigure}
\newcommand{\mailto}[1]{E-mail:~\href{mailto:#1}{\texttt{#1}}}

\usepackage[textwidth=.95\marginparwidth,
  backgroundcolor=lightgray,
  linecolor=orange,
  textsize=scriptsize]{todonotes}

\numberwithin{equation}{section}

\makeatletter
\@mparswitchfalse  % Place todo notes in wide margin.
\def\listtodoname{Todo List}
\def\listoftodos{%
  \section*{\listtodoname}\mbox{~}\par
  \@starttoc{tdo}
  }
\makeatother

\usepackage{natbib}
\bibpunct{(}{)}{;}{a}{,}{,}
\newcommand*{\doi}[1]{DOI: \href{http://dx.doi.org/\detokenize{#1}}{\detokenize{#1}}}

\makeatletter
\renewcommand\tableofcontents{%
  \section*{\contentsname}\mbox{~}\par
    \@starttoc{toc}
    }
\newcommand*\l@section[2]{%
  \ifnum \c@tocdepth >\z@
    \addpenalty\@secpenalty
    \setlength\@tempdima{1.5em}%
    \begingroup
      \parindent \z@ \rightskip \@pnumwidth
      \parfillskip -\@pnumwidth
      \leavevmode \bfseries
      \advance\leftskip\@tempdima
      \hskip -\leftskip
      #1\nobreak\hfil \nobreak\hb@xt@\@pnumwidth{\hss #2}\par
    \endgroup
  \fi}
\newcommand*\l@subsection{\@dottedtocline{2}{1.5em}{2.3em}}
\newcommand*\l@subsubsection{\@dottedtocline{3}{3.8em}{3.2em}}
\newcommand*\l@paragraph{\@dottedtocline{4}{7.0em}{4.1em}}
\newcommand*\l@subparagraph{\@dottedtocline{5}{10em}{5em}}
\makeatother

\newsiamthm{remark}{Remark}

\usepackage{dsfont}
\newcommand{\R}{\mathds{R}}
\newcommand{\minim}{\mathop{\mathrm{minimize}}}
\newcommand{\minimize}[1]{\displaystyle\minim_{#1}}
\newcommand{\maxim}{\mathop{\mathrm{maximize}}}

\newcommand{\st}{\mathop{\mathrm{subject\ to}}}
\newcommand{\Null}{\mathop{\mathrm{Null}}}

\usepackage{boxedminipage}
\newenvironment{btheorem}{%
   \setlength{\fboxsep}{4pt}%
   \par\hbox{}\noindent%
   \begin{boxedminipage}{\textwidth}
     \begin{theorem}}
    {\end{theorem}
    \end{boxedminipage}
    }
\newenvironment{blemma}{%
   \setlength{\fboxsep}{4pt}%
   \par\hbox{}\noindent%
   \begin{boxedminipage}{\textwidth}
     \begin{lemma}}
    {\end{lemma}
    \end{boxedminipage}
    }

\usepackage{algorithmicx}
\usepackage{algpseudocode}
\usepackage{enumitem}
\usepackage{mathtools}
\usepackage{bbm}
\usepackage{subcaption}
\usepackage{xspace}
\usepackage{tabularx}
\usepackage{booktabs}
\usepackage{subcaption}
\usepackage{standalone}

\crefname{subsection}{section}{sections}

\title{Implementing a smooth exact penalty function
    \\ for general constrained nonlinear optimization%
       \thanks{Version of \today}}
\author{%
  Ron Estrin%
  \thanks{Institute for Computational and Mathematical Engineering, Stanford  University, Stanford, CA 94305-4042 (\mailto{restrin@stanford.edu}).}
  \and
  Michael P.~Friedlander%
    \thanks{Department of Computer Science,
    University of British Columbia, Vancouver V6T 1Z4, BC, Canada
    (\mailto{mpf@cs.ubc.ca}). The work of this author was supported by
    ONR award N00014-17-1-2009.}
  \and
  Dominique Orban%
    \thanks{GERAD and Department of Mathematics and Industrial Engineering,
    \'Ecole Polytechnique, Montr\'eal, QC, Canada
    (\mailto{dominique.orban@gerad.ca}). The work of this author was
    supported by NSERC Discovery Grant 299010-04.}
  \and
  \hbox{Michael A. Saunders}%
    \thanks{Department of Management Science and Engineering, Stanford University, Stanford, CA 94305-4121 (\mailto{saunders@stanford.edu}). Research partially supported by the National Institute of General Medical Sciences of the National Institutes of Health [award U01GM102098].}}
\input{./macros}

\newcommand{\phis}{\phi_{\sigma}}

\newcommand{\ys}{y_{\sigma}}
\newcommand{\Ys}{Y_{\sigma}}
\newcommand{\gLag}{g\subl}
\newcommand{\hLag}{H\subl}
\newcommand{\Hs}{H_{\sigma}}
\newcommand{\gs}{g_{\sigma}}
\newcommand{\ws}{w_{\sigma}}
\newcommand{\Ws}{W_{\sigma}}
\newcommand{\gsy}{\gs^y}

\newcommand{\Ss}{S_{\sigma}}
\newcommand{\Rs}{R_{\sigma}}

\newcommand{\critcone}{\Cscr}

\newcommand{\augmats}{\bmat{I & Q^{1/2} A \\ A^T Q^{1/2}}}
\newcommand{\augmatn}{\bmat{I & A \\ A^T Q}}
\newcommand{\augmatnT}{\bmat{I & QA \\ A^T}}

\newcommand*{\LNLQ}{\hbox{\small LNLQ}\xspace}

\newcommand*{\CRAIG}{\hbox{\small CRAIG}\xspace}

\def\ind{\mathbbm{1}}

\newcommand{\mr}[1]{\multirow{2}{*}{#1}}

\begin{document}

\nolinenumbers
\maketitle

\begin{center} \it Dedicated to Roger Fletcher \end{center}

\thispagestyle{firstpage}
\pagestyle{myheadings}

\begin{abstract}
  We build upon \citet{EstrFrieOrbaSaun:2019a} to develop a general
  constrained nonlinear optimization algorithm based on a smooth
  penalty function proposed by
  \cite{Fletcher:1970,Fletcher:1973b}. Although Fletcher's approach has historically
  been considered impractical, we show that the computational kernels required
  are no more expensive than those in other widely accepted methods for nonlinear
  optimization. The main kernel for evaluating the penalty function and its
  derivatives solves structured linear systems.
  When the matrices are available
  explicitly, we store a single factorization each iteration.  Otherwise, we
  obtain a factorization-free optimization algorithm by solving each linear system
  iteratively.  The penalty function shows promise in cases where the linear systems
  can be solved efficiently, e.g., PDE-constrained optimization problems when
  efficient preconditioners exist.  We demonstrate the merits of the approach, and
  give numerical results on several PDE-constrained and standard test problems.
\end{abstract}

\section{Introduction}
\label{sec:introduction}
We consider a penalty-function approach for solving general constrained nonlinear optimization problems
\begin{equation}
\label{eq:nlp}
\tag{NP}
\begin{aligned}
  \minimize{x \in \R^n} \qquad & f(x) &&
\\ \st \qquad & c(x)=0 &&:\ y
\\ & \ell \le x \le u &&:\ z,
\end{aligned}
\end{equation}
where $f: \R^n \to \R$ and $c: \R^n \to \R^m$ are smooth
functions $(m \le n)$, the $n$-vectors $\ell$ and $u$ provide
(possibly infinite) bounds on $x$, and
$y \in \R^m$, $z \in \R^n$ are Lagrange multipliers associated with the equality constraints and bounds respectively.
\cite{EstrFrieOrbaSaun:2019a} describe factorization-based and factorization-free implementations of a smooth exact penalty method proposed by \cite{Fletcher:1970} to treat
equality constraints. Here, we generalize our implementation to
problems with both equality and bound constraints, and hence to problems with general inequality constraints.

% A smooth exact
% penalty function is used to eliminate the constraints $c(x)=0$ but
% leave the bound constraints explicit. This approach builds upon the
% method for equality-constrained problems described by

\citeauthor{Fletcher:1970}'s penalty function for equality constraints is the Lagrangian
\begin{equation}
    \label{eq:def-L}
    L(x,y) = f(x) - y\T c(x),
\end{equation}
in which the vector $y = \ys(x)$ is treated as a function of $x$ dependent on a parameter $\sigma > 0$.
\citet{Fletcher:1973b} proposes an extension to inequality constraints that exhibits nonsmoothness when constraint activities change.
The penalty function \eqref{eq:def-L} was long considered too costly for practical use \citep{Bert:1975, ConnGoulToin:2000, NocedalW:2006}, and the nonsmooth extension to inequality constraints further impacted its practicality. 

% This penalty function is minimized over just the bound
% constraints. The penalty function for equality constrained problems
% was proposed by \cite{Fletcher:1970}, with  in . 

We demonstrate that a certain smooth extension of \citeauthor{Fletcher:1970}'s penalty function yields a practical implementation for inequality-constrained optimization, by showing that the computational kernels are no more expensive than those in other widely accepted methods for nonlinear optimization, such as sequential quadratic programming.

The extended penalty function is {\em exact} because KKT points
of~\eqref{eq:nlp} are KKT points of the penalty problem for all
values of $\sigma$ larger than a finite threshold $\sigma^*$.  The
main computational kernel for evaluating the penalty function and its
derivatives is the solution of certain structured linear systems. We show how to solve the systems efficiently by factorizing a single matrix each iteration (if the matrix is available explicitly) and reusing the factors to evaluate the penalty function and its derivatives.
We also provide a {\em factorization-free} implementation in which linear systems are solved iteratively.
This makes the penalty function particularly applicable to certain problem classes such as PDE-constrained problems, where excellent preconditioners exist (e.g., those based on \citep{ReesDollWath:2010,StolWath:2012,Ridz:2013}); see \cref{sec:numerical-experiments}.

The advantage of smooth exact penalty functions is that they lead to conceptually simpler algorithms compared to traditional methods for constrained problems.
The original problem is replaced by a single smooth bound-constrained problem with a sufficiently large penalty parameter.
This avoids complicated heuristics to trade-off primal and dual feasibility, and can avoid the need for primal feasibility restoration stages or composite-step methods. Further, because our penalty is smooth and we can compute a sufficiently accurate Hessian approximation, second-order methods with fast local convergence may be used.

\subsection*{Paper outline}
\label{sec:outline}

We follow the structure of \cite{EstrFrieOrbaSaun:2019a}.
We introduce the penalty function in \cref{sec:penalty-function}, and discuss its relationship with existing approaches in \cref{sec:related-work}.  We give the penalty function's properties and derive
an explicit threshold for the penalty parameter in
\cref{sec:prop-penalty-funct}.  In \cref{sec:comp-penalty-funct} we
show how to evaluate the penalty function and its derivatives
efficiently.  We discuss an extension to maintain linear constraints in \cref{sec:explicit-constraints}.  Practical
considerations pertaining to the penalty function appear in
\cref{sec:practical-considerations}. We apply the penalty approach to standard and PDE-constrained problems
in \cref{sec:numerical-experiments}, and discuss future research directions in
\cref{sec:discussion}.

\section{The proposed penalty function}
\label{sec:penalty-function}

For \eqref{eq:nlp}, we propose the penalty function
\begin{equation}
  \label{eq:phi}
  \phis(x) := f(x) - c(x)\T \ys(x) = L(x, \ys(x)),
\end{equation}
where $\ys(x)$ are Lagrange multiplier estimates defined with other items as
\begin{align}
   \label{eq:3}
   \ys(x) &:= \textstyle{\argmin_{y}}\
                   \half\norm{A(x)y - g(x)}_{Q(x)}^2 + \sigma c(x)\T y,
                   & g(x) &:= \nabla f(x),
\\ \label{eq:8}
   A(x)        &:= \nabla c(x) = \bmat{g_1(x) \ \cdots\ g_m(x)},
                   & g_i(x) &:= \nabla c_i(x),
\\ \label{eq:Y}
   \Ys(x) &:= \nabla \ys(x).
\end{align}
Note that $A$ and $\Ys$ are $n$-by-$m$ matrices. We define an $n$-by-$n$
diagonal matrix $Q(x) = \textrm{diag}(q_i(x_i))$ 
with $\omega \in \R_+^n$, $\omega < u - \ell$,
and
\begin{equation}
\begin{aligned}
\label{eq:q-def}
{
    q_i(x_i) := \begin{cases} 1            & \mbox{if $\ell_i = -\infty$ and $u_i = \infty$}, \\
                            %   x_i - \ell_i & \mbox{if } \ell_i > -\infty,\ u_i = \infty, \\
                            %   u_i - x_i    & \mbox{if } \ell_i = -\infty,\ u_i < \infty, \\
                              \half (u_i - \ell_i) - \tfrac{1}{4}\omega_i - \tfrac{1}{4\omega_i} \left( 2x_i - u_i - \ell_i\right)^2 & \mbox{if } |u_i + \ell_i - 2x_i| \le \omega_i, \\
                              \min\{x_i - \ell_i, u_i - x_i\} & \mbox{otherwise.}\end{cases}
}
\end{aligned}
\end{equation}%\smarttodo{These conditions are mutually exclusive}
The diagonal of $Q(x)$ is a smooth approximation of $\min \{x-\ell, u - x\}$, and $\omega$ 
controls the smoothness. We use $\omega_i = \min\{ 1, \half (u_i - \ell_i) \}$.  Note that $Q(x)$ is nonnegative on $[\ell, u]$. We describe this function in more detail below.

We assume that~\eqref{eq:nlp} satisfies the following conditions:

% \smallskip

\begin{enumerate}[label=(A\arabic*)]
    \item \label{assump:c3} $f$ and $c$ are $\Cscr_3$.
    \item \label{assump:licq} The linear independence constraint
      qualification (LICQ) is satisfied for stationary points
      and all $x$ satisfying $\ell < x < u$.  LICQ is satisfied at $x$ if
      $$
      \left\{\nabla c_i(x), e_j \mid x_j \in \{\ell_j , u_j \}, \ i
        \in [m], \ j \in [n] \right\}$$ is linearly independent, where
      $e_j$ is the $j$th column of the identity matrix, and
      $[n] := \{1, 2, \dots, n\}.$
    \item \label{assump:strict-complementarity} Stationary points
      satisfy strict complementarity.  If $(\xstar, \ystar, \zstar)$
      is a stationary point, exactly one of $\zstar_j$ and
      $\min\{\xstar_j-\ell_j,u_j-\xstar_j\}$ is zero
      for all $j \in [n]$.
    \item The problem is feasible.  That is, there exists $x$ such that $\ell \le x \le u$ and $c(x) = 0$, with $\ell_j < u_j$ for all $j \in [n]$. We assume fixed variables have been eliminated from the problem.
    %\smarttodo{Do you mean no fixed vars? (Fixed vars aren't treated as bounds, they would be treated as equality constraints.)}
\end{enumerate}

\smallskip Assumption~\ref{assump:c3} ensures that $\phis$ has two
continuous derivatives and is typical for smooth exact penalty
functions \citep[Proposition 4.16]{Bert:1982}. However, at most two
derivatives of $f$ and $c$ are required to implement this penalty
function in practice (see section~\ref{sec:comp-second-deriv}).
Assumption~\ref{assump:licq} guarantees that $\Ys(x)$ and $\ys(x)$ are
uniquely defined; \ref{assump:strict-complementarity} provides
additional regularity to ensure that the threshold penalty parameter
$\sigma^*$ is well defined.

The basis of our approach is to solve
\begin{equation}
\label{eq:penalty-problem}
\tag{PP}
\minimize{x \in \R^n} \quad \phis(x) \quad
\st \quad \ell \le x \le u :\ z
\end{equation}
instead of \eqref{eq:nlp}.
We purposely set $z$ to be the Lagrange multiplier for the bound constraints of both \eqref{eq:nlp} and \eqref{eq:penalty-problem} because, as we show, they are equal at a solution.

\subsection{The scaling matrix}
\label{sec:q-description}

\begin{figure}[ht]
    \centering
    \includestandalone[mode=buildnew,page=1,scale=1]{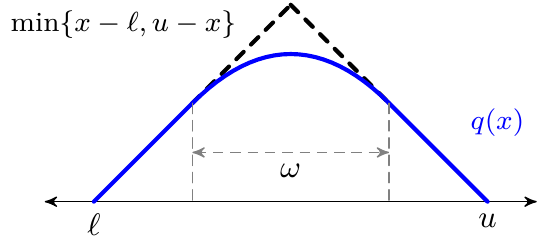}
    \caption{Plot of $q(x)$, a smooth approximation of $\min \{ x - \ell, u - x \}$.}
    \label{fig:q}
\end{figure}

The diagonal entries of the scaling matrix $Q(x)$ are smooth
approximations of the complementarity function $\min\{x-\ell, u-x\}$ \citep{Chen:2000}.
% , who gives a survey of several smoothing
% functions.
\cref{fig:q} plots $q(x)$ with finite $\ell$ and $u$.

The definition of $\ys(x)$ \eqref{eq:3} can be interpreted as a smooth
approximation of the complementarity conditions in the first-order KKT
conditions \eqref{eq:dual-complement}--\eqref{eq:dual-ineq2}
below. The role of $Q(x)$ is therefore to ensure that the partial derivatives of
the Lagrangian corresponding to indices of inactive bounds are zero.
Similar smoothing strategies can be found in the complementarity constraint
literature \citep{Anit:2000,Leyf:2006}.
%Several researchers have provided approaches for dealing with
%complementarity constraints in a similar fashion by introducing
%\smarttodo{Is it related? We don't have complementarity constraints. (We are approximating the complementarity conditions in the KKT conditions using the definition of $\ys(x)$)}
%nonlinear constraints---see for example \cite{Anit:2000} and
%\cite{Leyf:2006}.

For $x \in \R$, the derivative of $q(x)$ is
%\smarttodo{Need to explain the ``else ifs''. Need else ifs in (2.5)? (The conditions in (2.5) are mutually exclusive so ``else ifs" unnecessary. Do you think the SISC crowd will be confused?)}
\begin{equation}
  \label{eq:q-prime}
  q'(x) =
  \begin{cases} 0 & \mbox{if $\ell = -\infty$ and $u = \infty$,} \\
                              -\tfrac{1}{\omega} \left( 2x + u - \ell \right) & \mbox{if $|u + \ell - 2x| \le \omega$,} \\
                              1 & \mbox{if $x - \ell < u - x$,} \\
                              {-1} & \mbox{if $x - \ell > u - x$.}
                            \end{cases}
\end{equation}

Note that the cases in \eqref{eq:q-prime} are not mutually exclusive, and should be checked top to bottom until a case is satisfied.
The choice of $q(x)$ is not unique because any smooth concave function
that is zero at $x_j\in \{\ell_j, u_j\}$ works in our framework. For
instance, if \(u_j - \ell_j\) is large, we could use a smooth approximation of
$\min \{ x_j - \ell_j, u_j - x_j, 1 \}$ to avoid
numerical issues that can arise if $x$ is far from its bounds.

\subsection{Notation}
\label{sec:notation}

Denote $\xstar$ as a local stationary point
%\smarttodo{Minimizer or stationary point?}
of \eqref{eq:nlp}, with corresponding dual solutions $\ystar$ and $\zstar$. At $\xstar$, define the set of active bounds as
%\smarttodo{$\Ascrbar$ not used later?}
\begin{equation}
\label{eq:active-set}
\Ascr(\xstar) := \{j \mid x_j \in \{\ell_j, u_j\}\}, %\qquad \Ascrbar(\xstar) := \{j \mid \ell_j < x_j < u_j\},
\end{equation}
and define the critical cones $C_{\phi}(\xstar, \zstar)$ and $C(\xstar, \zstar)$ as
\begin{subequations}
\label{eq:critical-cone}
\begin{align}
\critcone_{\phi}(\xstar, \zstar) &:= \left\{ p\, \middle\vert \begin{array}{cl} p_j = 0 & \mbox{if } \zstar_j \ne 0 \\ p_j \ge 0 & \mbox{if } \xstar_j = \ell_j \\ p_j \le 0 & \mbox{if } \xstar_j = u_j\end{array} \right\},
\\ \critcone(\xstar, \zstar) &:= \left\{ p \in \critcone_{\phi}(\xstar,\zstar) \,\middle\vert\, A(\xstar)\T p = 0 \right\}.
\end{align}
\end{subequations}
Observe that by \ref{assump:strict-complementarity}, $\critcone_{\phi}(\xstar, \zstar) = \left\{ p \mid p_j = 0 \mbox{ if } \zstar_j \ne 0 \right\}$, so $p \in \critcone_{\phi}(\xstar,\zstar)$ if and only if $p = Q(\xstar)^{1/2} \pbar$ for some $\pbar \in \R^n$.
%\smarttodo{I don't understand this, could you explain more? z has disappeared. If $x_j = u_j$, $q_j(x) = 0$. How does that imply $p_j \leq 0$? (By strict complementarity (A3), if $x_j$ is on a bound, then $z_j \ne 0$, so $p_j = 0$. Also $q_j(x_j) = 0$, so $(Q(x)\pbar)_j = 0$. If $x_j$ is interior, then $q_j(x_j) \ne 0$, so we can find $\pbar_j$ such that $q_j(x_j) \pbar_j = p_j$.)}

Let $g(x) = \nabla f(x)$, $H(x)=\nabla^2 f(x)$, $g_i(x) = \nabla c_i(x)$, $H_i(x)=\nabla^2 c_i(x)$, and define
%\begin{subequations} \label{eq:grad-hess-sig}
%\begin{align}
%  g_\sigma(x) & := g(x) - A(x) \ys(x), \label{eq:grad-sig},
%  \qquad
%  H_\sigma(x) & := H(x) - \sum_{i=1}^m [\ys(x)]_i H_i(x),\label{eq:hess-sig}
%\end{align}
%\end{subequations}
\begin{equation} \label{eq:grad-hess-sig}
\begin{aligned}
    \gLag(x,y) &:= g(x) - A(x) y,
    \qquad&
    \gs(x) &:= \gLag(x,\ys(x)),
\\  \hLag(x,y) &:= H(x) - \sum_{i=1}^m y_i H_i(x),
    \qquad&
    \Hs(x) &:= \hLag(x,\ys(x))
\end{aligned}
\end{equation}
as the gradient and Hessian of $L$ at $(x,y)$ or $(x,\ys(x))$.
We define the matrix operators
%\smarttodo{What are the indices $i$ in $R(x,v)$? Should they be $1, \dots, n$? (What do you mean?)}
%\smarttodo{Define $g_i$ in $S$ and $T$ (Defined in (2.3))}
\begingroup
\allowdisplaybreaks
\begin{align*}
   R(x,v) &:= \nabla_x[Q(x)v] = \nabla_x \bmat{q_1(x_1) v_1 \\ \vdots \\ q_n(x_n) v_n} = \diag \left( \bmat{ q_1^{\prime}(x_1) v_1 \\ \vdots \\ q_n^{\prime}(x_n) v_n} \right),
\\ S(x,v) &:= \nabla_x[A(x)\T v]
            = \nabla_x \bmat{g_1(x)^T v \\ \vdots \\ g_m(x)^T v}
            =          \bmat{v^T H_1(x) \\ \vdots \\ v^T H_m(x)},
\\
T(x,w) &:= \nabla_x[A(x) w]
            = \nabla_x \left[ \sum_{i=1}^m w_i g_i(x) \right]
            =                 \sum_{i=1}^m w_i H_i(x),
\end{align*}
\endgroup
where $v \in \R^n$, $w \in \R^m$, and $T$ is a symmetric matrix. The operation of multiplying the adjoint of $S$ with a vector $w$ is described by
\begin{align*}
%   S(x,v)   w &= \bmat{v^T H_1(x) w \\ \vdots \\ v^T H_m(x) w}, \\
S(x,v)^T w &= \left[\sum_{i=1}^m w_i H_i(x)\right] v = T(x,w) v = T(x,w)^T v \,.
%\\T(x,w)^T v &= \left[\sum_{i=1}^m w_i H_i(x)\right] v = T(x,w) v \,.
\end{align*}
%with $v \in \R^n$, $w \in \R^m$.
If $A_Q(x) = Q(x)^{1/2}A(x)$ has full rank $m$, the operators
\begin{equation}
  \label{eq:proj}
  P(x) := A_Q(x)\big(A_Q(x)\T A_Q(x)\big)\inv A_Q(x)^T
  \text{and}
  \Pbar(x) := I - P(x)
\end{equation}
define orthogonal projectors onto $\range(A_Q(x))$ and its complement
respectively. More generally, for a matrix $M$, we define $P_M$ and $\Pbar_M$ as the orthogonal projectors onto $\range(M)$ and $\Null(M)$
%\smarttodo{Usually, we write $\Null(M)$!? (I use $\ker(M)$ for the space, and $\Null(M)$ for the dimension; not insisting though)}
respectively. % (which may or may not be a function of $x$).
%If $M$ has full column-rank, $M^{\dagger} = (M\T M)^{-1} M^T$ denotes the Moore-Penrose pseudoinverse of $M$.
%\smarttodo{No more daggers later in the text?}

Unless otherwise indicated, $\|\cdot\|$ is the 2-norm for vectors and matrices. For $M$ positive definite, $\|u\|^2_M = u\T M u$ is the energy-norm. For square matrices $M$, define $\lambda_{\max}(M)$ as its largest eigenvalue. Define $\ind$ as the vector of all ones of size dictated by the context.

\section{Related work on penalty functions for inequality constraints}
\label{sec:related-work}

Penalty functions have long been used to solve constrained problems by
replacing constraints with functions that penalize infeasibility.
% We provide a brief overview of other smooth exact penalty methods and their relation to \eqref{eq:penalty-problem}.
\citet[\S 1.1]{EstrFrieOrbaSaun:2019a} give an overview of other smooth exact penalty methods for equality constrained optimization and their relation to \eqref{eq:penalty-problem}.
% , which can also be applied in this case with minor modification, so we do not repeat it here.
A more detailed overview is given by \citet{PilloGrippo:1984}, \citet{ConnGoulToin:2000}, and \citet{NocedalW:2006}.

When $\ell = 0$ and $u = \infty$, \cite{Fletcher:1973b} proposes
the penalty function
\begin{align*}
    \psi_\sigma(x) &:= f(x) - c(x)\T \ys(x) - z_{\sigma}(x)\T x,
\\  \ys(x) ,z_{\sigma}(x) &:= \argmin_{\{y \in \R^m,\, z \ge 0\}} \half \|A(x) y + z - g(x)\|_2^2 + \sigma c(x)\T y
\end{align*}
and minimizes $\psi_\sigma$ unconstrained.
Although $\psi_\sigma$ is exact and continuous, it is nonsmooth because of the bound
constraints on $z$: active-set changes on those
bounds correspond to non-differentiable points for $\psi_\sigma$.
Solving the penalty problem requires a method for nonsmooth problems, and \cite{Mara:1978}
observes that nonsmooth merit functions may result in slow convergence.

Since \cite{Fletcher:1973}, there has been significant work on smooth exact penalty methods that handle inequality constraints \citep{PilloGrippo:1984,PillGrip:1985,BoggTollKear:1991,AnitZava:2014}. Many approaches replace the inequality constraints with equalities using squared slacks \citep{Bert:1982}, at which point the equality constrained problem is solved via a smooth exact penalty approach. (This is one approach for deriving $\phis$ and \eqref{eq:3}; however, it is also possible to derive it directly from the first-order KKT conditions.) The penalty function in these cases is the augmented Lagrangian, which either keeps the dual variables explicit and penalizes the gradient of the Lagrangian \citep{AnitZava:2014}, or expresses the dual variables as a function of $x$ \citep{PilloGrippo:1984}. Our penalty function \eqref{eq:phi} takes the latter approach but defines this parametrization differently from previous approaches; rather than introducing additional dual variables for the bounds in \eqref{eq:3}, we change the norm of the least-squares problem according to the distance from the bounds, to approximate the complementarity conditions of first-order KKT points.

\section{Properties of the penalty function}
\label{sec:prop-penalty-funct}
%\smarttodo{I removed the subsection headers to avoid very short sections}

In this section, we show how $\phis(x)$ naturally expresses
the optimality conditions of~\eqref{eq:nlp}. We also give explicit
expressions for the threshold value of the penalty parameter $\sigma$.

% \subsection{Derivatives of the penalty function}
% \label{sec:derivatives-of-phi}

As in \citep{EstrFrieOrbaSaun:2019a}, the gradient and Hessian of $\phis$ may be written as
\begin{subequations} \label{eq:phi-grad-hess}
  \begin{align}
    \nabla  \phis(x)  & = g_\sigma(x) - \Ys(x) c(x),
    \label{eq:phi-grad}
    \\
    \nabla^2 \phis(x) & = H_\sigma(x)- A(x) \Ys(x)^T
                                    - \Ys(x) A(x)^T
                                    - \nabla_x \left[\Ys(x) c \right],
    \label{eq:phi-hess}
  \end{align}
\end{subequations}
where %$H(x)=\nabla^2 f(x)$ and $H_i(x)=\nabla^2 c_i(x)$ are the
%Hessians of the objective and each constraint function, respectively.
the last term $\nabla_x[\Ys(x) c]$ %in the expression for $\nabla^2\phis$
purposely drops the argument on $c$ to emphasize that this gradient is
made on the product $\Ys(x) c$ with $c:=c(x)$ held fixed.  This
term involves third derivatives of $f$ and $c$, and as we shall see,
it is convenient and computationally efficient to ignore it.  We leave
it unexpanded.

% \subsection{Optimality conditions}
% \label{sec:optim-cond}

The penalty function $\phis$ is closely related to the (partial) Lagrangian~\eqref{eq:def-L}.
To make this connection clear, we define the Karush-Kuhn-Tucker (KKT)
optimality conditions for~\eqref{eq:nlp} in terms of those of~\eqref{eq:penalty-problem}.
%
%The KKT conditions are usually defined via the Lagrangian; e.g., see \cite[Ch.~12]{NocedalW:2006}.
From the definition of $\phis$ and $\ys$ and~\eqref{eq:phi-grad-hess},
% definitions of $\phis$ and its gradient and Hessian
% the following definitions are equivalent
% to the KKT conditions.
we have the following definition.

\begin{definition}[First-order KKT points of \eqref{eq:nlp}]
\label{def:kkt-1}
  The point $(\xstar, \zstar)$ is a first-order KKT point of~\eqref{eq:nlp} if for any $\sigma \geq 0$ the following hold:
  %\smarttodo{Does ``for any'' mean $\forall$ or $\exists$? ($\exists$ because if it holds for some value of $\sigma$, it will hold for all of them.)}
  \begingroup
  \allowdisplaybreaks
  \begin{subequations} \label{eq:1st-order}
  \begin{align}
    \label{eq:6}
    \ell \le \xstar &\le u,
  \\c(\xstar)                 &= 0, \label{eq:primal-feas}
  \\\nabla\phis(\xstar) &= \zstar, \label{eq:dual-feas-phi}
  \\\zstar_j &= 0, \qquad \mbox{if } j \notin \Ascr(\xstar), \label{eq:dual-complement}
  \\\zstar_j &\ge 0, \qquad \mbox{if } \xstar_j = \ell_j, \label{eq:dual-ineq1}
  \\\zstar_j &\le 0, \qquad \mbox{if } \xstar_j = u_j. \label{eq:dual-ineq2}
  \end{align}
  \end{subequations}
  \endgroup
  Then $\ystar:=\ys(\xstar)$ is the
  Lagrange multiplier of~\eqref{eq:nlp} associated with $\xstar$. Note that by \ref{assump:strict-complementarity}, inequalities \eqref{eq:dual-ineq1} and \eqref{eq:dual-ineq2} are strict.
  %\smarttodo{Should we make them strict then? (I don't think so, because the definition is general. The rest is a consequence of our assumptions)}
\end{definition}

\begin{remark}
If \eqref{eq:1st-order} holds for some $\sigma \ge 0$, it necessarily holds for all $\sigma \ge 0$ because $c(\xstar) = 0$.
Also, the point $(\xstar, \zstar)$ is a first-order KKT point of~\eqref{eq:penalty-problem} if for any $\sigma \geq 0$, \eqref{eq:6} and \eqref{eq:dual-feas-phi}--\eqref{eq:dual-ineq2} hold.
\end{remark}

%We can similarly derive second-order optimality conditions based on
%the Hessian of $\phis$.

\begin{definition}[Second-order KKT point of \eqref{eq:nlp}]
\label{def:kkt-2}
  The first-order KKT point $(\xstar,\zstar)$ satisfies the second-order necessary KKT condition
  for~\eqref{eq:nlp} if for any $\sigma \geq 0$,
  \begin{equation} \label{eq:2nd-order-nec}
  \begin{aligned}
    p^T \nabla^2\phis(\xstar) p \ge 0
    \quad
    \hbox{for all $p \in \critcone(\xstar, \zstar)$. }
  \end{aligned}
  \end{equation}
%   \begin{equation} \label{eq:2nd-order-suf}
%     p^T \nabla^2\phis(\xstar) p > 0
%     \quad
%     \hbox{for all $p\ne0$ such that }
%     A(\xstar)\T p = 0.
%   \end{equation}

%   The first-order KKT point $(\xstar,\zstar)$ satisfies the second-order necessary KKT condition
%   for~\eqref{eq:penalty-problem} if for any $\sigma \geq 0$,  \smarttodo{I merged the remark with the definitions. Ok? (I think we should keep them separate. Might get confusing because they're not part of the definition)}
%   \begin{equation} \label{eq:2nd-order-necphi}
%   \begin{aligned}
%     p^T \nabla^2\phis(\xstar) p \ge 0
%     \quad
%     \hbox{for all $p \in \critcone_{\phi}(\xstar,\zstar)$. }
%   \end{aligned}
%   \end{equation}
%   The conditions~\eqref{eq:2nd-order-nec} and~\eqref{eq:2nd-order-necphi} are sufficient if the inequality is strict.
Condition~\eqref{eq:2nd-order-nec} is sufficient if the inequality is strict.
\end{definition}

\begin{remark}
  %If \eqref{eq:primal-feas} is omitted, \cref{def:kkt-1} corresponds
  %to first-order KKT points of \eqref{eq:penalty-problem}.
  %Similarly,
  If $(\xstar,\zstar)$ is a first-order KKT point for \eqref{eq:penalty-problem}, then replacing $\critcone(\xstar,\zstar)$ by
  $\critcone_{\phi}(\xstar,\zstar)$ in \cref{def:kkt-2} corresponds to
  second-order KKT points of \eqref{eq:penalty-problem}.
\end{remark}

The second-order KKT condition says that at a second-order KKT point
of \eqref{eq:penalty-problem}, $\phis$ has nonnegative curvature
along directions in the critical cone
$\critcone_\phi(\xstar, \zstar)$. We now show that at $\xstar$,
increasing $\sigma$ increases curvature only along the normal cone
to the equality constraints. We derive a threshold value for $\sigma$ beyond which
that $\phis$ has nonnegative curvature even when
$A(\xstar)\T p \ne 0$, as well as a condition on $\sigma$ that ensures that
stationary points of \eqref{eq:penalty-problem} are primal
feasible. For a given first- or second-order KKT triple
$(\xstar,\ystar,\zstar)$ of \eqref{eq:nlp}, we define
\begin{equation}
  \label{eq:14}
  \sigma^* := \half\lambda^+_{\max}\left( P(\xstar) Q(\xstar)^{1/2} \hLag(\xstar,\ystar) Q(\xstar)^{1/2} P(\xstar) \right),
\end{equation}
where $\lambda^+_{\max}(\cdot) = \max\{\lambda_{\max}(\cdot), 0\}$.
The following lemmas are similar to those of \cite{EstrFrieOrbaSaun:2019a}. Indeed, if the bounds are absent then $Q(x) = I$ and we recover the same results as in \cite{EstrFrieOrbaSaun:2019a}.
% \smarttodo{I think without strict complementarity,
%   we get a similar threshold, but it's nastier. Or maybe not.}
%\smarttodo{If the lemma and theorem are variants of something in Fletcher1, we should say it. (I agree; do you think that works?)}

\begin{blemma}
If $c(x) \in \range \left(A(x)\T Q(x) \right)$, then $\ys(x)$ satisfies
\begin{equation}
    \label{eq:lin-sys-ys}
    A(x)\T Q(x) A(x) \ys(x) = A(x)\T Q(x) g(x) - \sigma c(x).
\end{equation}
Furthermore, if $Q(x) A(x)$ has full rank, then
\begin{equation}
    \label{eq:AAYTx}
    \begin{aligned}
    &A(x)\T Q(x) A(x) \Ys(x)\T
    \\ &\qquad = A(x)\T \left[ Q(x)\Hs(x) - \sigma I + R(x, \gs(x)) \right] + S(x,Q(x)\gs(x)).
    \end{aligned}
\end{equation}
\end{blemma}

\begin{proof}
  For any $x$, the necessary and sufficient optimality conditions
  for~\eqref{eq:3} give~\eqref{eq:lin-sys-ys}. For brevity, let
  everything be evaluated at the same point $x$ and drop the argument
  $x$ from all operators. By differentiating both sides
  of~\eqref{eq:lin-sys-ys}, we obtain
\[
   S(Q A \ys) + A\T \left[ R(A \ys) + Q T(\ys) + Q A \Ys^T \right]
 = S(Q g) + A\T \left[R(g) + Q H - \sigma I \right].
\]
The derivative exists because $\ys(x)$ is well-defined in a neighbourhood of $x$ if $Q(x)A(x)$ is full-rank.
By rearranging the above and using definitions \eqref{eq:grad-hess-sig}, we obtain \eqref{eq:AAYTx}.
\end{proof}
% \begin{remark}
% By \ref{assump:licq}, $c(x) \in \range\left( A(x)\T Q(x) \right)$ is always satisfied.
% \end{remark}

\begin{btheorem}[Threshold penalty value]
  \label{thm:threshold}
  Suppose $(\xbar,\zbar)$ is a first-order KKT point for
  \eqref{eq:penalty-problem} with $Q(\xbar)^{1/2}A(\xbar)$ full-rank, and let $(\xstar,\ystar,\zstar)$ be a
  second-order necessary KKT point for \eqref{eq:nlp}.  Then
  %\vspace*{-12pt}
  \begin{subequations}
  \begin{align}
    \label{eq:17}
    \sigma > \norm{A(\xbar)^T Q(\xbar) \Ys(\xbar)}
       &\quad\Longrightarrow\quad c(\xbar) = 0; % \cref{def:kkt-1};
% \\  \label{eq:13}
%     \sigma \ge \norm{A(\xstar_1) \Ys(\xstar_1)^T}
%       &\quad\Longrightarrow\quad \sigma \ge \sigma^*;
\\  \label{eq:16}
    p^T \nabla^2\phis(\xstar) p \geq 0 \text{for all} p \in \critcone_{\phi}(\xstar,\zstar)
       &\quad\!\Longleftrightarrow\quad \sigma \ge \sigmabar,
  \end{align}
  \end{subequations}
  where $\sigmabar = \half\lambda_{\max}\left( P(\xstar) Q(\xstar)^{1/2} \hLag(\xstar,\ystar) Q(\xstar)^{1/2} P(\xstar) \right)$ is defined in~\eqref{eq:14}.
  The consequence of~\eqref{eq:17} is that \(\xbar\) is a first-order KKT point for~\eqref{eq:nlp}.
  If $\xstar$ is second-order sufficient, the inequalities in~\eqref{eq:16} hold strictly.
  Observe that $\sigma^* = \max\{\sigmabar,0\}$.
\end{btheorem}

\begin{proof}

Proof of~\eqref{eq:17}:
% Because $\xbar$ is a first-order KKT point for \eqref{eq:penalty-problem}, we need only show that $c(\xbar) = 0$.
By~\eqref{eq:dual-feas-phi}--\eqref{eq:dual-ineq2}, $Q(\xbar) \nabla \phis(\xbar) = 0$, so that
\[
  Q(\xbar) g(\xbar) = Q(\xbar) A(\xbar) \ys(\xbar) + Q(\xbar) \Ys(\xbar) c(\xbar).
\]
Substituting~\eqref{eq:lin-sys-ys} evaluated at $\xbar$ into this equation yields, after
simplifying,
%\smarttodo{Need to assume $c \in R(A^TQ)$? If that wasn't the case, $\Ys$ wouldn't exist. This also follows from (A2)}
\[
  A(\xbar)\T Q(\xbar) \Ys(\xbar) c(\xbar) = \sigma c(\xbar).
\]
Taking norms of both sides and using the triangle inequality gives the
inequality $\sigma\norm{c(\xbar)}\le\norm{A(\xbar)^T Q(\xbar)
  \Ys(\xbar)}\,\norm{c(\xbar)}$, which implies that
$c(\xbar) = 0$.

Proof of~\eqref{eq:16}: Because $\xstar$ satisfies first-order
conditions~\eqref{eq:1st-order}, we have $\ystar = \ys(\xstar)$ and
$Q(\xstar) g_\sigma(\xstar)=0$, independently of $\sigma$.
Therefore $S(\xstar,Q(\xstar)\gs(\xstar)) = 0$. We drop the
arguments from operators that take $x$ as input and assume that they
are all evaluated at $\xstar$. 
%It follows from~\eqref{eq:AAYTx}, $\hLag(\xstar,\ystar) = \Hs$, and the definition of the projector $P := P(\xstar)$ that
By premultiplying \eqref{eq:AAYTx} by $(A_Q^\dagger)^T = Q^{1/2}A(A\T Q A)^{-1}$ and postmultiplying by $Q^{1/2}$, using $\hLag(\xstar,\ystar) = \Hs$, and the definition of $P := P(\xstar)$, we have
\begin{align}
  Q^{1/2}A \Ys^T Q^{1/2} &= (A_Q^\dagger)^T A^T(Q\hLag(\xstar, \ystar)Q^{1/2}-\sigma I + R(\gs))Q^{1/2}
\\&= PQ^{1/2} \hLag(\xstar, \ystar)Q^{1/2} - \sigma P + (A_Q^\dagger)^T A R(\gs)Q^{1/2}   \label{eq:11}
\end{align}
Observe that if $p\in \critcone_{\phi}(\xstar,\zstar)$, then $p = Q^{1/2} \pbar$ for some $\pbar \in \critcone_{\phi}(\xstar,\zstar)$.
%\smarttodo{$\bar{p} \in \R^n$? (Yes, because Q is square)}
Because $Q^{1/2} \gs = 0$, we have $R(\gs) Q^{1/2} = 0$.
Therefore using \eqref{eq:phi-hess}, \eqref{eq:11}, $c(\xstar) = 0$, and the relation $P + \Pbar = I$, we have
\begin{align*}
   p^T \nabla^2 \phis(\xstar) p \ge 0 &\Leftrightarrow \pbar^T Q^{1/2} \left( \Hs - A \Ys\T - \Ys A^T \right) Q^{1/2} \pbar \ge 0
\\   &\Leftrightarrow \pbar^T \left( Q^{1/2} \Hs Q^{1/2} - PQ^{1/2} \Hs Q^{1/2} - Q^{1/2} \Hs Q^{1/2}P + 2\sigma P \right) \pbar
\\   &\Leftrightarrow \pbar^T \left( \Pbar Q^{1/2} \Hs Q^{1/2} \Pbar - P Q^{1/2} \Hs Q^{1/2} P + 2\sigma P \right) \pbar \ge 0.
\end{align*}
Now, because $\Pbar \pbar \in \Null (A^T Q^{1/2})$ implies that $Q^{1/2} \Pbar \pbar \in \critcone(\xstar, \zstar)$, the first term above is nonnegative according to \cref{def:kkt-2}. It follows that $\sigma$ must be sufficiently large that $2 \sigma P - P Q^{1/2} \hLag(\xstar,\ystar) Q^{1/2} P \succeq 0$, which is equivalent to $\sigma \geq \sigmabar$.
\end{proof}

As in \citet[Theorem 4]{EstrFrieOrbaSaun:2019a}, \eqref{eq:16} shows
that if $\xstar$ is a second-order KKT point of~\eqref{eq:nlp}, there
exists a threshold value $\sigmabar$ beyond which $\xstar$ is also a
second-order KKT point of~\eqref{eq:penalty-problem}. As penalty parameters are typically nonnegative, we treat $\sigma^* = \max\{\sigmabar,0\}$ as the threshold. Note that this
result does not preclude the possibility that there exist minimizers
of the penalty function---for any value of $\sigma$---that are not
minimizers of \eqref{eq:nlp}. However, these are rarely encountered in
practice. Further, we can add a quadratic penalty term that, under
certain conditions, ensures that KKT points of~\eqref{eq:penalty-problem} are feasible for~\eqref{eq:nlp}
\citep[\S 3.3]{EstrFrieOrbaSaun:2019a}.

\section{Evaluating the penalty function}
\label{sec:comp-penalty-funct}

The main challenge in evaluating $\phis$ and its gradient is the
solution of the shifted weighted-least-squares problem~\eqref{eq:3}
needed to compute $\ys(x)$, and computation of the gradient
$\Ys(x)$. We show below that it is possible to compute matrix-vector
products $\Ys(x)v$ and $\Ys(x)\T u$ by solving structured linear
systems involving the same matrix. We show that this linear system may
be either symmetric or unsymmetric, and discuss the tradeoffs
between both approaches. In either case, if direct methods are to be
used, only a single factorization that defines the
solution~\eqref{eq:3} is required for all products.

For this section, it is convenient to drop the arguments on various functions and assume they are all evaluated at a point $x$ for some parameter $\sigma$.  For example,
%\smarttodo{You already use this notation in the proof of the lemma and theorem. Why not introduce it before the lemma?}
\(
\ys = \ys(x), \ A = A(x), \ \Ys = \Ys(x), \ H_\sigma = H_\sigma(x), \ S_\sigma = S(x,Q(x)g_\sigma(x)), \Rs = R(x,\gs(x))\), etc.  We express~\eqref{eq:AAYTx} using the shorthand notation
\begin{equation}
  \label{eq:ATAY}
  A\T Q A \Ys^T = A^T (Q \Hs - \sigma I + \Rs) + \Ss.
\end{equation}
We first describe how to compute products $\Ys u$ and $\Ys\T v$, then how to put those pieces together to evaluate the penalty function and its derivatives.

Every quantity of interest can be computed by solving a symmetric or
unsymmetric linear system and combining the solution with the
derivatives of the problem data. Typically it is preferable to solve
symmetric systems; however, we find
that additional Jacobian products are then needed. The additional cost may be negligible, but this matter
becomes application-dependent. We therefore present both options,
beginning with the symmetric case.

There are many ways to construct the right-hand sides of the linear
systems presented below. One consideration is that inversions with the
diagonal matrix $Q^{1/2}$ should be avoided---even though the diagonal
of $Q$ will be assumed strictly positive because of the use of an
interior method (see \cref{sec:practical-considerations}), numerical
difficulties may arise near the boundary of the feasible set if
$Q^{1/2}$ contains small entries and is
inverted. %, which can impede convergence.

\subsection{Computing \texorpdfstring{$\mathbf{\Ys u}$}{Yu}}

It follows from \eqref{eq:ATAY} that for a given $m$-vector $u$,
\begin{align*}
    \Ys u &= (\Hs Q - \sigma I + \Rs) A (A\T Q A)^{-1} u + \Ss\T (A\T Q A)^{-1} u.
% \\     &= \Hs Q^{1/2}v + (\sigma I - \Rs) Aw - \Ss\T w
% \\     &= (\Hs - \sigma I + \Rs) \vbar - \Ss\T w,
\end{align*}

%\subsubsection{Symmetric linear system}
Let $w = -(A^T Q A)^{-1} u$ and $v = -Q^{1/2} Aw$, so that $v$ and $w$
are the solution of the symmetric linear system
\begin{equation}
  \label{eq:aug-Yu-sym}
  \augmats \bmat{v\\w} = \bmat{
0\\u}.
\end{equation}
Then $\Ys u = \Hs Q^{1/2}v + (\sigma I - \Rs) Aw - \Ss\T w$. \Cref{alg:Yu-sym} formalizes this process.
\begin{algorithm}[hbt]
  \caption{%
    Computing the matrix-vector product $\Ys u$
    \label{alg:Yu-sym}
  }
  \begin{algorithmic}[1]
    \State $(v,w) \gets \hbox{solution of~\eqref{eq:aug-Yu-sym}}$
    \State \Return $\Hs Q^{1/2}v + (\sigma I - \Rs) Aw - \Ss\T w$
  \end{algorithmic}
\end{algorithm}

% \subsubsection{Unsymmetric linear system}
% Define $w = -(A^T Q A)^{-1} u$, and $\vbar = -Aw$, so that $\Ys u =  (\Hs - \sigma I + \Rs) \vbar - \Ss\T w$,
% and $\vbar$ and $w$ are the solutions of the unsymmetric linear system
% \begin{equation}
%   \label{eq:aug-Yu-nsym}
%   \augmatn \bmat{\vbar \\ w} = \bmat{
% 0\\u}
% \end{equation}
% \autoref{alg:Yu-nsym} formalizes the process.

% \begin{algorithm}[hbt]
%   \caption{%
%     Computing the matrix-vector product $\Ys u$ (unsymmetric system)
%     \label{alg:Yu-nsym}
%   }
%   \begin{algorithmic}[1]
%     \State $(\vbar,w) \gets \hbox{solution of~\eqref{eq:aug-Yu-nsym}}$
%     \State \Return $(\Hs - \sigma I + \Rs) \vbar - \Ss\T w$
%   \end{algorithmic}
% \end{algorithm}

\subsection{Computing \texorpdfstring{$\mathbf{\Ys\T v}$}{Y'v}}
Again from \eqref{eq:ATAY}, multiplying both sides by
$v$ gives
%\smarttodo{Is it $+ \Rs$?}
\begin{align*}
    \Ys\T v &= (A\T Q A)^{-1} A^T (Q H_\sigma-\sigma I + \Rs)v + (A\T Q A)^{-1} S_\sigma v.
\end{align*}
%\subsubsection{Symmetric linear system}
The product $u=\Ys\T v$ is part of the solution of the system
\begin{equation}
  \label{eq:aug-YTv-sym}
  \augmats \bmat{r \\ u} =
  \bmat{Q^{1/2} \Hs v \\ A\T(\sigma I - \Rs)v -S_\sigma v}.
\end{equation}
\cref{alg:YTv-sym} formalizes the process.

\begin{algorithm}[hbt]
  \caption{%
    Computing the matrix-vector product $\Ys\T v$
    \label{alg:YTv-sym}
  }
  \begin{algorithmic}[1]
    \State Evaluate $Q^{1/2} \Hs v$ and $A\T(\sigma I + \Rs)v -S_\sigma v$
    \State $(r,u) \gets \hbox{solution of~\eqref{eq:aug-YTv-sym}}$
    \State \Return $u$
  \end{algorithmic}
\end{algorithm}

% \subsubsection{Unsymmetric linear system}
% The product $u=\Ys\T v$ is the solution of
% \begin{equation}
%   \label{eq:aug-YTv-nsym}
%   \augmatnT \bmat{\rbar \\ u} =
%   \bmat{(Q H_\sigma-\sigma I - \Rs)v \\ -S_\sigma v}.
% \end{equation}
% \cref{alg:YTv-nsym} formalizes the process.

% \begin{algorithm}[hbt]
%   \caption{%
%     Computing the matrix-vector product $\Ys\T v$
%     \label{alg:YTv-nsym}
%   }
%   \begin{algorithmic}[1]
%     \State Evaluate $(Q H_\sigma-\sigma I - \Rs)v$ and $S_\sigma v$
%     \State $(r,u) \gets \hbox{solution of~\eqref{eq:aug-YTv-nsym}}$
%     \State \Return $u$
%   \end{algorithmic}
% \end{algorithm}

\subsection{Unsymmetric linear system}

We briefly comment on how to use unsymmetric systems in place of \eqref{eq:aug-Yu-sym} and \eqref{eq:aug-YTv-sym}. We can compute products of the form $\Ys u = (\Hs - \sigma I + \Rs) \vbar - \Ss\T w$ (where $w = -(A^T Q A)^{-1} u$ and $\vbar = -Aw$), and products $u = \Ys\T v$ by solving the respective linear systems:
%\smarttodo{Another $+ \Rs$?}
\begin{equation}
  \label{eq:aug-nsym}
  \augmatn \bmat{\vbar \\ w} = \bmat{
    0\\u} \text{and}
  \augmatnT \bmat{\rbar \\ u} =
  \bmat{(Q H_\sigma-\sigma I - \Rs)v \\ -S_\sigma v}.
\end{equation}
Algorithms \ref{alg:Yu-sym} and \ref{alg:YTv-sym} can then be appropriately modified to use the above linear systems.

\subsection{Computing multipliers and first derivatives}

The multiplier estimates $\ys$ and Lagrangian gradient can be obtained from one of the following linear systems:
\begin{equation}
    \label{eq:aug-mult}
    \augmats \bmat{d \\ \ys} = \bmat{Q^{1/2}g \\ \sigma c} \qquad \mbox{or} \qquad \augmatn \bmat{\gs \\ \ys} = \bmat{g \\ \sigma c}.
\end{equation}
Observe that in the unsymmetric case we obtain $\gs$ immediately.
The symmetric system yields $d = Q^{1/2} \gs$. As noted earlier, computing $\gs \gets Q^{-1/2} d$ may amplify errors when the diagonal entries of $Q$ are approaching zero.
An alternative would be to compute $\gs \gets g - A\ys$, which costs an extra Jacobian product.

The penalty gradient $\nabla \phis = \gs - \Ys c$ can then be computed using $\gs$ and computing $\Ys c$ via \cref{alg:Yu-sym} or its unsymmetric variant.% or \cref{alg:Yu-nsym}.

\subsection{Computing second derivatives} \label{sec:comp-second-deriv}

We approximate $\nabla^2 \phis$ from \eqref{eq:phi-hess} using the
same approaches as \citet{EstrFrieOrbaSaun:2019a}:
\begin{subequations}
\label{eq:hess-approx}
\begin{align}
    \nabla^2\phis &\approx B_1 := \Hs - A\Ys^T - \Ys A^T \label{eq:hess-approx-1}
\\  &\phantom{\approx B_1 :}= \Hs - \widetilde{P} (Q \Hs + \Rs - \sigma I)  - (\Hs Q + \Rs - \sigma I) \widetilde{P} \nonumber
\\  &\phantom{\approx B_1 :}\qquad - A(A\T Q A)^{-1} S_{\sigma} -  S_{\sigma}^T (A\T Q A)^{-1}A \nonumber
\\  &\approx B_2 := \Hs - \widetilde{P} (Q \Hs + \Rs - \sigma I)  - (\Hs Q + \Rs - \sigma I) \widetilde{P}, \label{eq:hess-approx-2}
\end{align}
\end{subequations}
where $\widetilde{P} = A(A\T Q A)^{-1}A$. The first approximation
drops the third derivative term $\nabla [\Ys c]$ in
\eqref{eq:phi-hess}, while the second approximation drops the term
$\Ss(x, Q \gs)$, because those terms are zero at a
solution. Thus, $B_1$ and $B_2$ can be interpreted as %a form of
Gauss-Newton approximations of $\nabla^2 \phis$. Using similar
arguments to those made by \citet[Theorem 2]{Fletcher:1973}, we expect
those approximations to result in quadratic convergence when
$f, c \in \Cscr_3$, and at least superlinear convergence when
$f,c \in \Cscr_2$.

Computing products with $B_1$ only requires products with $\Ys$ and $\Ys^T$, which can be handled by Algorithms \ref{alg:Yu-sym} and \ref{alg:YTv-sym}.
To compute a product  $\widetilde{P} u$, we can solve
\begin{equation}
    \augmats \bmat{p \\ q} = \bmat{0 \\ A\T u}
    \enspace \mbox{or} \enspace
    \augmatn \bmat{\pbar \\ q} = \bmat{u \\ 0},
    \qquad \widetilde{P} u = -Aq.
\end{equation}
As before, using the unsymmetric system avoids an additional Jacobian
product, which may be negligible compared to solving an unsymmetric
system.

\subsection{Solving the augmented linear system}

We comment on various approaches for solving the necessary linear systems
\begin{equation}
\label{eq:aug-generic}
    \Kscr \bmat{p \\ q} = \bmat{w \\ z}, \quad \mbox{ where } \quad \Kscr = \augmats \mbox{ or } \augmatn.
\end{equation}
This is the most computationally intensive step in our approach. %, and care should be taken when choosing what method to apply (and to decide whether to solve the symmetric or unsymmetric system).
%It should be noted that when using direct methods, a single factorization can be computed per iteration and used for every system solve
Note that with direct methods, a single factorization is needed to evaluate $\phis$ and its derivatives.

\citet[\S 4.5]{EstrFrieOrbaSaun:2019a} describe several approaches for solving the symmetric system (using both direct and iterative methods), so we do not repeat this discussion here. For unsymmetric systems, %if direct methods are to be used,
any sparse factorization of $\Kscr$ may be used;
also, we could factorize $Q^{1/2}A$ with a Q-less QR factorization and use the (refined) semi-normal equations \citep{PaigBjor:1994} as in the symmetric case (as long as multiplications with $Q^{-1/2}$ are avoided).

If iterative methods are used, the unsymmetric system requires unsymmetric iterative methods such as GMRES \citep{SaadSchu:1986}, SPMR \citep{EstrGrei:2018}, or QMR \citep{FreuNach:1991}, where the choice of method depends on considerations such as short- vs.\ long-recurrence, available preconditioners, or robustness. Note that preconditioners approximating $\Pscr \approx A\T Q A$ apply to both the symmetric and unsymmetric systems; however, unsymmetric %systems have additional flexibility such as the ability to use inexact preconditioner solves, which short-recurrence symmetric methods may not be able to handle.
solvers may allow inexact preconditioner solves, while short-recurrence symmetric solvers may not.

If optimization solvers that accept inexact function and derivative
evaluations are used (e.g., \citet[\S 8--9]{ConnGoulToin:2000} or
\cite{HeinRidz:2014}), the results of \citet[\S
7]{EstrFrieOrbaSaun:2019a} apply here as well; that is, bounding the
residual norm of the linear systems is sufficient to bound the
function and derivative evaluation error up to a constant (under mild
assumptions). This is useful in cases where solving the linear
system exactly every iteration is prohibitively expensive. Further,
when the symmetric system is used, it is possible to use methods that
upper bound the solution error. % (e.g., \CRAIG
%\smarttodo{Arioli didn't invent CRAIG. (The citation was meant for his error-bounding work on CRAIG)}
% \citep{Ario:2013} or \LNLQ
% \citep{EstOrbSaun:2018}) when an underestimate of the smallest singular
% value of the preconditioned Jacobian is available.
For example, \citet{Ario:2013} develop error bounds for \CRAIG \citep{Craig:1955}, and \cite{EstOrbSaun:2019LNLQ} develop error bounds for \LNLQ when an underestimate of the smallest singular
value of the preconditioned Jacobian is available.

\section{Maintaining explicit constraints}
\label{sec:explicit-constraints}
We consider a variation of \eqref{eq:nlp} where some of the constraints $c(x)$ are easy to maintain explicitly; for example, linear equality constraints. %Many solvers can efficiently handle linear constraints \citep{X,XX,XXX} \smarttodo{get examples or drop}.
We show below that maintaining subsets of constraints explicitly decreases the threshold penalty parameter $\sigma^*$ in \eqref{eq:14}.
% We can then maintain feasibility for this subset of constraints, the contours of $\phis$ are simplified, and as we show soon, the threshold penalty parameter $\sigma^*$ is decreased.
%We discuss the case where some of the constraints are linear, but it is possible to extend the theory to any type of constraint.
Instead of~\eqref{eq:nlp}, consider the problem with explicit linear equality constraints
\begin{equation}
    \label{eq:nlp-exp}
    \tag{NP-EXP}
    \minimize{x \in \R^n} \enspace f(x) \enspace\st\enspace c(x)=0, \enspace B\T x = d, \enspace \ell \le x \le u,
\end{equation}
where $c(x) \in \R^{m_1}$ and $B\T x = d$ with $B \in \R^{n \times m_2}$, so that $m_1 + m_2 = m$.
We assume that \eqref{eq:nlp-exp} at least satisfies \ref{assump:licq}, so that $B$ has full column rank.
We define the penalty problem as
\begin{equation}
\label{eq:penalty-problem-exp}
\begin{aligned}
    \minimize{x \in \R^n} \enspace \phis (x) &:= f(x) - c(x)^T \ys(x) \enspace\st\enspace B\T x = d, \enspace \ell \le x \le u,
\\  \bmat{\ys(x) \\ \ws(x)} &:= \argmin_{y,w} \half \norm{A(x) y + Bw - g(x)}_{Q(x)}^2 + \sigma \bmat{c(x) \\ B\T x - d}\T \bmat{y \\ w},
\end{aligned}
\end{equation}
which is similar to \eqref{eq:penalty-problem} except that the linear constraints are not penalized in $\phis(x)$, and the linear constraints are explicitly present. Another possibility is to penalize the linear constraints as well, while keeping them explicit; however, this introduces additional nonlinearity in $\phis$.
Further, if all constraints are linear, it is desirable for the penalty function to reduce to~\eqref{eq:nlp-exp}.

For a given first- or second-order KKT solution $(\xstar,\ystar)$, the threshold penalty parameter becomes
\begin{align}
    \sigma^* :\!\!&= \half \lambda^+_{\max}\left( \Pbar_{Q^{1/2} B} P_{Q^{1/2}C} Q^{1/2} \hLag(\xstar,\ystar) Q^{1/2} P_{Q^{1/2}C} \Pbar_{Q^{1/2} B} \right)     \label{eq:threshold-eq}
    \\ &\leq \half\lambda^+_{\max}\left( P_{Q^{1/2}C} Q^{1/2} \hLag(\xstar,\ystar) Q^{1/2} P_{Q^{1/2}C} \right),     \label{eq:non-threshold-eq}
 \end{align}
where $Q := Q(\xstar)$, $C := \bmat{A(x^*) & B}$ is the Jacobian for all constraints. %, and $\Bbar := Q^{1/2}B$. 
Inequality \eqref{eq:non-threshold-eq} holds because $\Pbar_{Q^{1/2} B}$ is an orthogonal projector. If the linear constraints were not explicit, the threshold value would be \eqref{eq:non-threshold-eq}. Intuitively, the threshold penalty value decreases 
%according to how much of the %top eigenspace
%space spanned by the eigenvectors corresponding to positive eigenvalues of the Lagrangian Hessian lies in the range of $Q^{1/2} B$,
%\smarttodo{What is the ``top eigenspace''? (Space spanned by positive eigenvectors. Is this clearer?}
%of the Lagrangian Hessian lies in the range of $B^T$, 
because positive semidefiniteness of $\nabla^2 \phis(\xstar)$ is only required on a lower-dimensional subspace.%along that space is guaranteed at a second-order point of \eqref{eq:penalty-problem-exp}.
%\smarttodo{Should ``by the underlying solver'' be ``near a second-order point''? (Yes, that works too. Is this better?}

The following result is analogous to \cref{thm:threshold}.% with the smaller threshold value.
%\smarttodo{What does it mean that $\xstar$ is a first and second-order point, respectively? Theorem~6 just says second-order point.}

\begin{btheorem}[Threshold penalty value with explicit constraints]
  \label{thm:threshold-eq}
  Suppose $(\xbar,\zbar)$ is a first-order necessary KKT point for \cref{eq:penalty-problem-exp},
  and let $(\xstar,\ystar,\zstar)$ be a second-order necessary KKT point for \eqref{eq:nlp-exp}.
  Define $\critcone_{\phi}^* := \critcone_{\phi}(\xstar,\zstar) \cap \Null(B\T)$, $Q := Q(\xbar)$, and $\Pbar_{Q^{1/2} B} := \Pbar_{Q^{1/2} B}(\xbar)$.
  Then %for all $p \in \critcone_{\phi}(\xstar,\zstar)$ such that $B\T p = 0$,
  \begin{subequations}
  \begin{align}
    \label{eq:34}
    \sigma > \norm{A(\xbar)^T Q^{1/2} \Pbar_{Q^{1/2} B} Q^{1/2} \Ys(\xbar)}
      &\quad\Longrightarrow\quad c(\xbar) = 0;
\\  \label{eq:36}
    p^T \nabla^2\phis(\xstar) p \succeq 0 \text{for all} p \in \critcone_{\phi}^*
      &\quad\!\Longleftrightarrow\quad \sigma \ge \sigmabar,
  \end{align}
  \end{subequations}
  where \(\sigmabar = \half \lambda_{\max}\left( \Pbar_{Q^{1/2} B} P_{Q^{1/2}C} Q^{1/2} \hLag(\xstar,\ystar) Q^{1/2} P_{Q^{1/2}C} \Pbar_{Q^{1/2} B} \right)\). Again, $\sigma^* = \max\{\sigmabar,0\}$.
  The consequence of~\eqref{eq:34} is that \(\xbar\) is a KKT point for~\eqref{eq:nlp}.
  If $\xstar$ is second-order sufficient, the inequalities
  in~\eqref{eq:36} hold strictly.
\end{btheorem}

The proof of the theorem, and details of evaluating the penalty function with explicit constraints, are given in \cref{app:explicit-constraints}. Although we only considered the linear case here, explicit nonlinear constraints can be handled with minor modifications.

\section{Practical considerations}
\label{sec:practical-considerations}

So far we have demonstrated that for sufficiently large $\sigma$, minimizers of \eqref{eq:nlp} are minimizers of \eqref{eq:penalty-problem},
%\smarttodo{What's the definition of ``weakly exact''? (Let's avoid that definition entirely. Just means that minimizers of penalty are not necessary minimizers of the original)}
%weakly exact smooth penalty function,
and we showed how to evaluate $\phis$ and its derivatives. By \ref{assump:licq} we know that $\phis$ is defined for all $\ell < x < u$. Although it may appear that any optimization solver can be applied to minimize \eqref{eq:penalty-problem}, the structure of $\phis$ lends itself more readily to certain types of solvers.

First, we recommend interior solvers rather than exterior or active-set methods. For $\phis(x)$ to be defined, we require that $Q(x) \succeq 0$ (thus disqualifying exterior point methods) and that $Q(x)^{1/2} A(x)$ have full column-rank (so that at most $n-m$ components of $x$ can be at one of their bounds). Even if \ref{assump:licq}
is satisfied, an active-set method may choose a poor active set that causes $\phis(x)$ to be undefined (or it may have too many active bounds). On the other hand, interior methods ensure that $Q(x) \succ 0$ and avoid this issue (at least until $x$ converges and approaches the bounds).

As in \citep{EstrFrieOrbaSaun:2019a}, Newton-CG type trust-region
solvers \citep{Stei:1983} should be used to solve
\eqref{eq:penalty-problem}. Products with approximations of
$\nabla^2 \phis(x)$ can be efficiently computed, but computing the
Hessian itself is not practical. Also, trust-region methods are better
equipped to deal with negative curvature than linesearch methods
($\phis$ typically has an indefinite Hessian). Finally, evaluating
$\phis$ at several points (such as during a linesearch) is expensive
because every evaluation requires solving a different linear
system. Given these considerations, a solver like KNITRO
\citep{ByrdNoceWalt:06} is ideal for solving
\eqref{eq:penalty-problem}.

It remains future work to determine a robust procedure for updating
$\sigma$ if it is too small (causing $\phis$ to be unbounded) or too
large (causing small steps to be taken). 
For the following experiments, we choose an initial $\sigma$ specific to each problem and keep it constant.
We also have the same heuristic available that is
discussed by \citet[\S 8]{EstrFrieOrbaSaun:2019a} to update $\sigma$,
which often works in practice.  %\smarttodo{Discuss updating $\sigma$ more?}

\section{Numerical experiments}
\label{sec:numerical-experiments}

We investigate the performance of Fletcher's penalty function on several PDE-constrained optimization problems and some standard test problems. For each test we use the stopping criterion
%\smarttodo{Notation clash with the operator T on p4 (Good catch! $N$?}
\begin{equation}
    \label{eq:experiment-stop}
    \begin{aligned}
      \|c(x)\|_\infty &\le \epsilon_p \\
      \|N(x) \gs(x)\|_\infty &\le \epsilon_d
    \end{aligned}
    \qquad\mbox{or}\qquad \|N(x) \nabla \phis(x)\|_\infty \le \epsilon_d,
\end{equation}
with $N(x) = \mbox{diag}(\min \{x-\ell, u-x, \ind \} )$, $\epsilon_p := \epsilon \left( 1 + \|x\|_{\infty} + \|c(x_0)\|_\infty \right)$, and $\epsilon_d := \epsilon \left( 1 + \|y\|_{\infty} + \|\gs(x_0)\|_\infty \right)$, where $x_0$ is the initial point, $y_0 = y_0(x_0)$, and $\epsilon = 10^{-8}$.

%Depending on the problem, the augmented systems \eqref{eq:aug-generic} are solved by either direct or iterative methods. For direct methods, we use the corrected semi-normal equations \citep[Appendix A.2.2]{EstrFrieOrbaSaun:2019a}.
For the standard test problems, we use the
%\smarttodo{Does that mean semi-normal eqns with one step of iterative refinement? (Yes; should we specify this?)}
semi-normal equations with one step of iterative refinement \citep{PaigBjor:1994}.
For the PDE-constrained problems, we use \LNLQ with the \CRAIG transfer point
%\smarttodo{Arioli didn't invent CRAIG (We use Arioli's bounds)}
\citep{EstOrbSaun:2019LNLQ, Craig:1955, Ario:2013} to solve the symmetric augmented system \eqref{eq:aug-generic} with preconditioner $\Pscr$ and two possible termination criteria:
\begingroup
\allowdisplaybreaks
\begin{subequations} \label{eq:system-solve-accuracy}
  \begin{align}
  \left\|\bmat{p^* \\ q^*} - \bmat{p^{(k)} \\ q^{(k)}} \right\|_{\widebar\Pscr^{\phantom{-1}}} &\le \eta \left\| \bmat{p^{(k)} \\ q^{(k)}}\right\|_{\widebar\Pscr}, \label{eq:terminate-error}
    \qquad \bar\Pscr := \bmat{I & \\ & \Pscr},
\\[3pt]  \left\| \Kscr \bmat{p^{(k)} \\ q^{(k)}} - \bmat{u \\ v} \right\|_{\widebar\Pscr^{-1}} &\le \eta \left\| \bmat{u \\ v}\right\|_{\widebar\Pscr^{-1}}, \label{eq:terminate-residual}
  \end{align}
\end{subequations}
\endgroup
which are based on the relative error and
the relative residual (obtained via \LNLQ \citep{EstOrbSaun:2019LNLQ}), respectively. We
can use \eqref{eq:terminate-error} when a lower bound on
$\sigma_{\min}(\Pscr^{-1/2} A)$ is available, which is the case in the
PDE-constrained optimization problems below.

We use KNITRO \citep{ByrdNoceWalt:06} to solve
\eqref{eq:penalty-problem}. For the PDE-constrained optimization problems, we set the penalty parameter to $\sigma = 10^t$, for the smallest $t$ that allowed KNITRO to converge. When $\phis$ is evaluated approximately
(for $\eta$ large), we use such solvers without modification, thus
pretending that the function and gradient are evaluated exactly.
The use of inexact linear solves is discussed in \cite[\S 7]{EstrFrieOrbaSaun:2019a};
the following experiments using inexactness are similar to those in \cite[\S 9]{EstrFrieOrbaSaun:2019a}.
%\smarttodo{What's the impact of this on the convergence? (Do you mean compared to properly handling inexactness? I suppose we technically don't guarantee convergence anymore...)}

\subsection{2D inverse Poisson problem}
\label{sec:inv-poisson}

\begin{figure}[t]
    \centering
    \begin{subfigure}{0.28\textwidth}
        \centering
        \includegraphics[height=103pt]{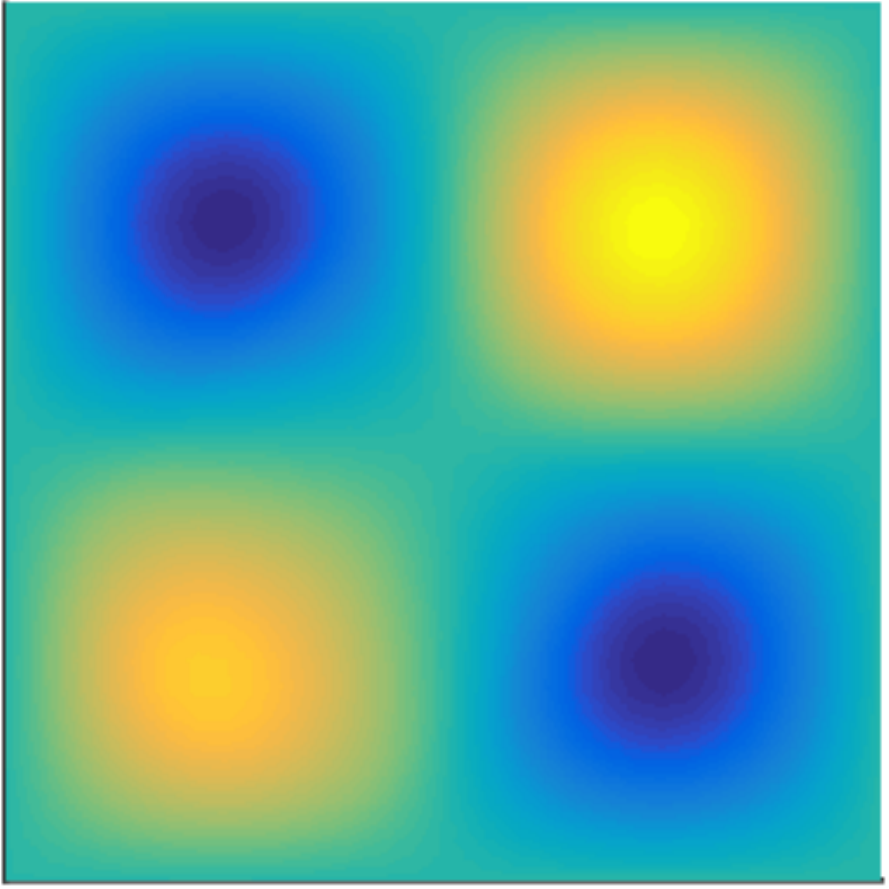}
        \caption{Target state $u_d$}
    \end{subfigure}
    ~
    \begin{subfigure}{0.28\textwidth}
        \centering
        \includegraphics[height=103pt]{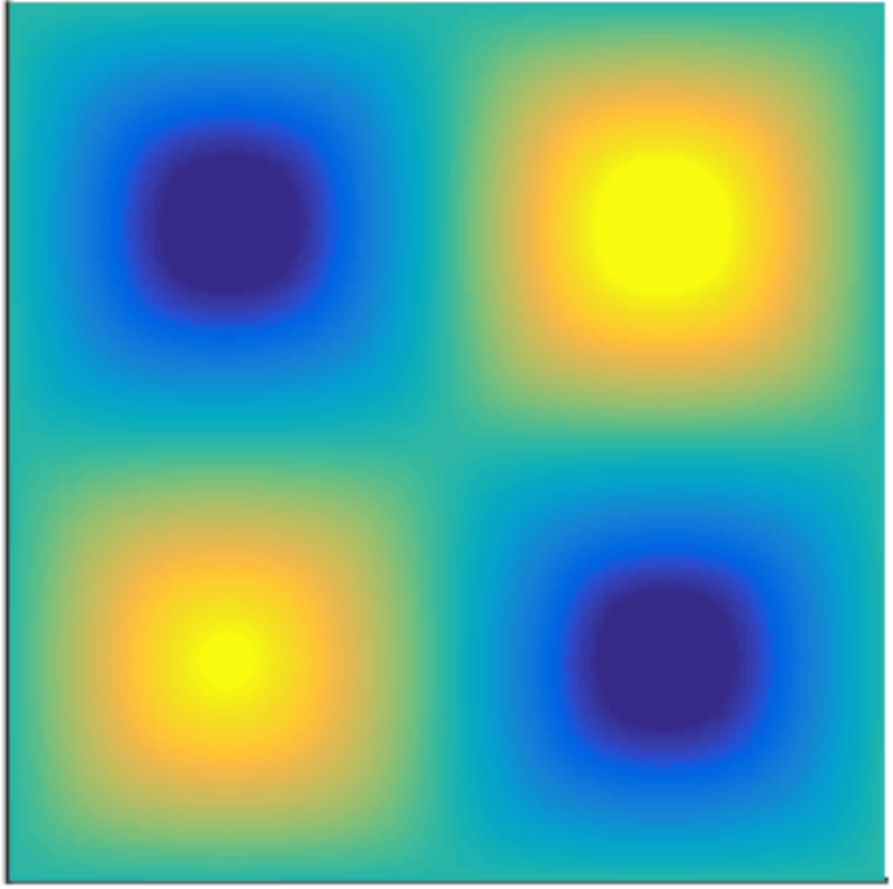}
        \caption{State (coarse mesh)}
    \end{subfigure}
    ~
    \begin{subfigure}{0.28\textwidth}
        \centering
        \includegraphics[height=103pt]{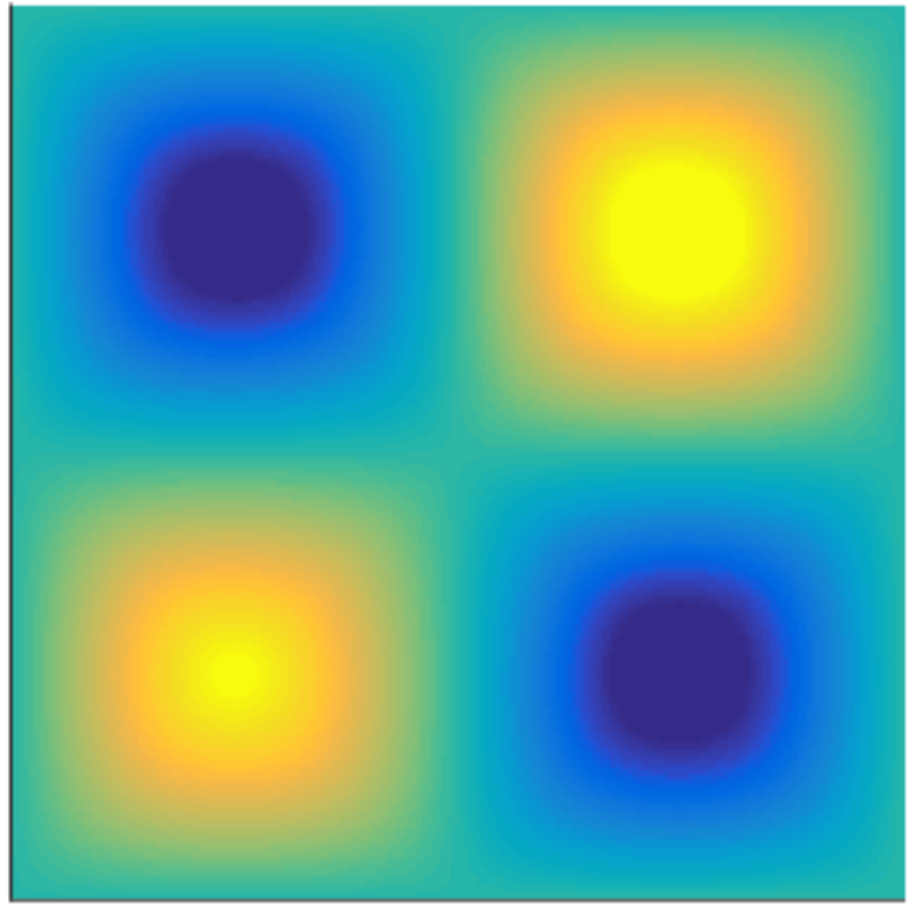}
        \caption{State (fine mesh)}
    \end{subfigure}
    ~
    \begin{subfigure}{0.08\textwidth}
        \centering
        \includegraphics[height=103pt]{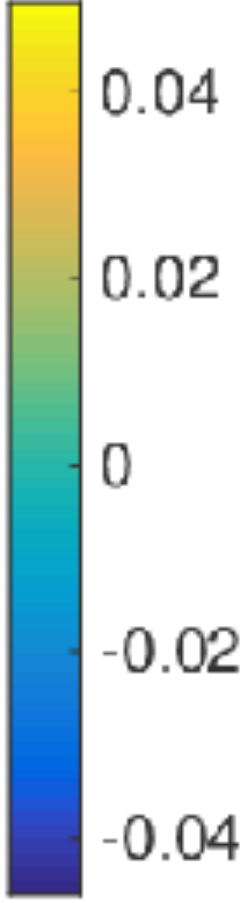}
        \caption*{}
    \end{subfigure}
    \\
    \begin{subfigure}{0.28\textwidth}
        \centering
        \includegraphics[height=103pt]{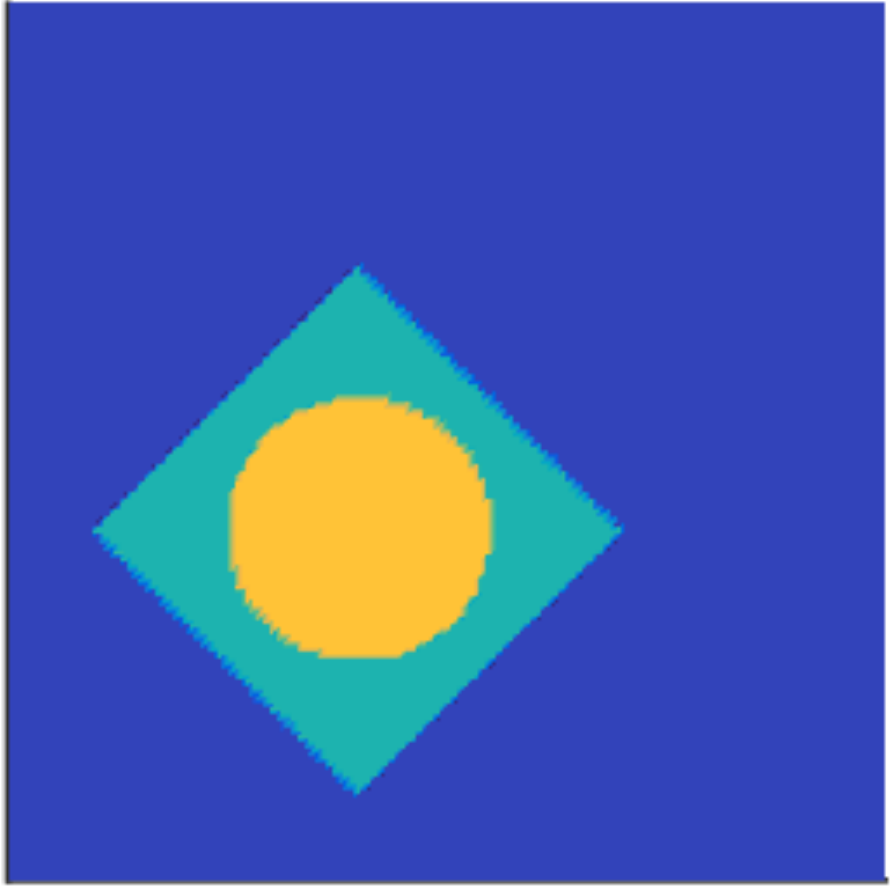}
        \caption{Target control $z_*$}
    \end{subfigure}
    ~
    \begin{subfigure}{0.28\textwidth}
        \centering
        \includegraphics[height=103pt]{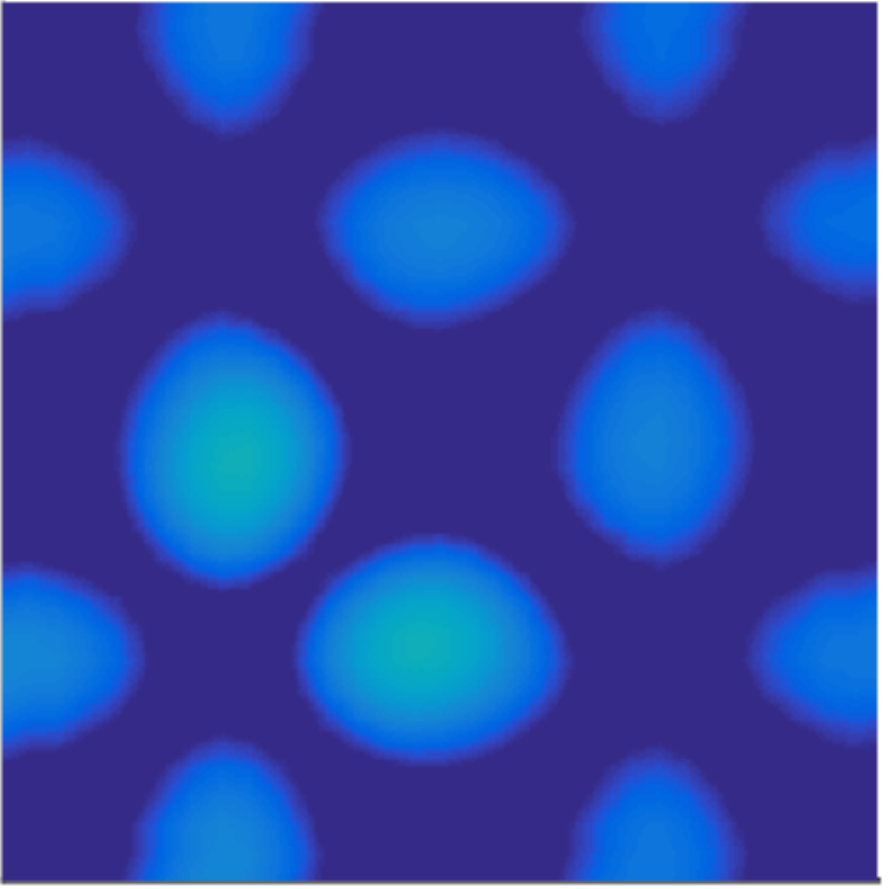}
        \caption{Control (coarse mesh)}
    \end{subfigure}
    ~
    \begin{subfigure}{0.28\textwidth}
        \centering
        \includegraphics[height=103pt]{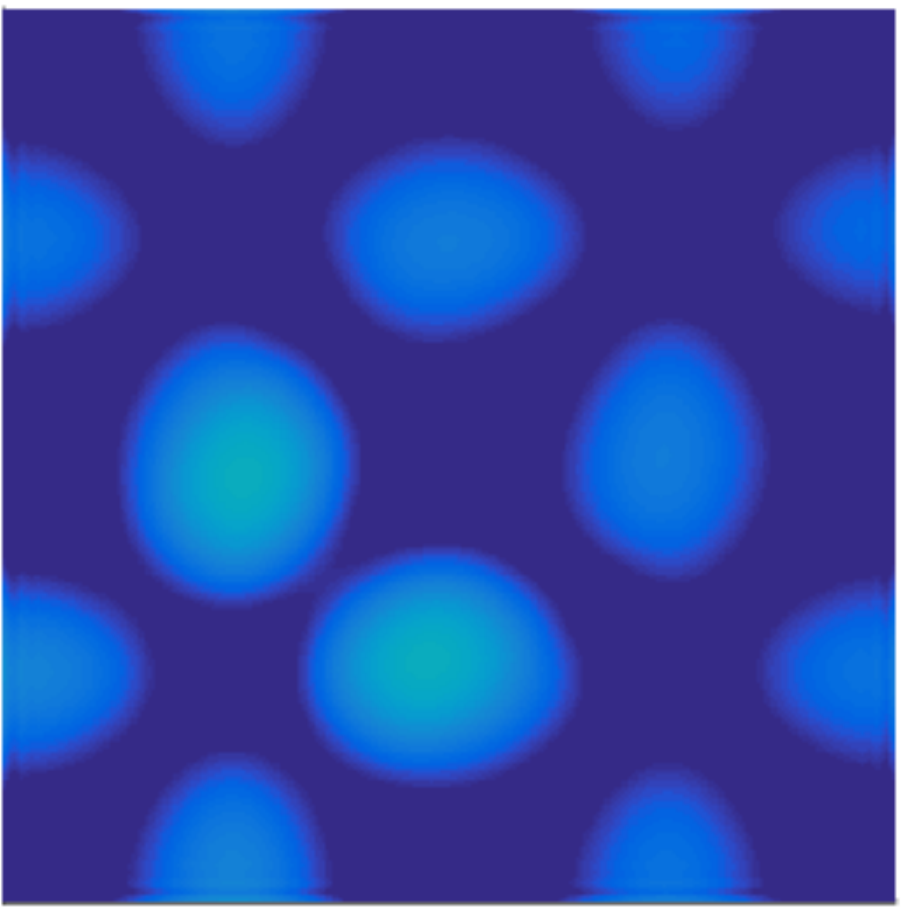}
        \caption{Control (fine mesh)}
    \end{subfigure}
    ~
    \begin{subfigure}{0.08\textwidth}
        \centering
        \hspace*{-10pt}
        \includegraphics[height=105pt]{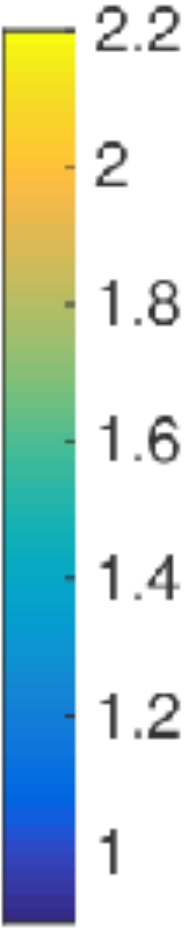}
        \caption*{}
    \end{subfigure}
    \caption{Target and computed states (top), and controls (bottom) for \eqref{eq:inverse-poisson}. Because the problem is ill-posed, the control is not exactly recovered, even though the state is well-matched.}
    \label{fig:inv_poisson}
\end{figure}
Let $\Omega = (-1,1)^2$ represent the physical domain and $H^1(\Omega)$ denote the Sobolev space of functions in $L^2(\Omega)$, whose weak derivatives are also in $L^2(\Omega)$. Let $H_0^1(\Omega) \subset H^1(\Omega)$ be the Hilbert space of functions whose value on the boundary $\partial \Omega$ is zero. We solve the following 2D PDE-constrained control problem:
\begin{equation}
\label{eq:inverse-poisson}
\begin{aligned}
    \minimize{u \in H_0^1(\Omega), \, z \in L^2(\Omega)} \enspace & \half \int_\Omega \left( u - u_d \right)^2 \mathrm{d}x + \tfrac{1}{2} \alpha \int_\Omega z^2 \, \mathrm{d}x
\\ \st \quad &
   \begin{array}[t]{rll}
     - \nabla \cdot (z \nabla u) = & \hspace{-.5em} h & \mbox{in } \Omega,
     \\ u = &\hspace{-.5em} 0 & \mbox{on } \partial \Omega,
     \\ z \ge & \hspace{-.5em} 0 & \mbox{in } \Omega.
   \end{array}
% \\  &\enspace - \nabla \cdot (z \nabla u) = h \qquad\mbox{in } \Omega,
% \\  &\enspace u = 0 \qquad\mbox{in } \partial \Omega.
\end{aligned}
\end{equation}
Let $c = (0.2,0.2)$ and define $S_1 = \{ x \mid \|x - c\|_2 \le 0.3 \}$ and $S_2 = \{x \mid \|x - c\|_1 \le 0.6 \}$. For a set $C$, define $I_C(x) = 1$ if $x \in C$ and 0 otherwise.
The target state $u_d$ is generated as the solution of the PDE with
$z_*(x) = 1 + 0.5 \cdot I_{S_1}(x) + 0.5\cdot I_{S_2}(x)$.

% Table 1
\begin{table}[t]
  \caption{Results from solving \eqref{eq:inverse-poisson} using KNITRO to solve \eqref{eq:penalty-problem} with various $\eta$ in \eqref{eq:terminate-error} (left) and \eqref{eq:terminate-residual} (right) to terminate the %augmented
  linear system solves. The top (resp.\ bottom) table records results for the smaller problem with $n=2050$, $m=1089$ (resp.\ larger problem with $n=20002$, $m=10201$). We record the number of function/gradient evaluations ($\#f,g$), Lagrangian Hessian (\#$Hv$), Jacobian (\#($Av$), and adjoint Jacobian (\#$A\T v$) products.}
  \label{tab:inv-poisson}
  \centering
  \begin{tabular}{@{} c |  c  c c  c  c | c c c c c @{}}
    \toprule
    $\eta$ & Its. & \#$f,g$ & \#$H v$ & \#$Av$ & \#$A\T v$ & Its. & \#$f,g$ & \#$H v$ & \#$Av$ & \#$A\T v$ \\
    \midrule
    $10^{-2\phantom{0}}$&46 &64 &2856 &8436 &8611 &67 &81 &4374 &12915 &13145 \\
    $10^{-4\phantom{0}}$&43 &55 &2168 &6642 &6796 &36 &51 &1458 &4642 &4781 \\
    $10^{-6\phantom{0}}$&35 &46 &2120 &6876 &7004 &29 &35 &1194 &4138 &4238 \\
    $10^{-8\phantom{0}}$&39 &50 &2322 &7833 &7973 &47 &71 &7062 &22150 &22340 \\
    $10^{-10}$          &37 &47 &2236 &8110 &8242 &43 &58 &3170 &11565 &11725 \\
    \midrule
    \multicolumn{11}{c}{} \vspace{-5pt}\\
    \toprule
    $10^{-2\phantom{0}}$&144 &176 &3662 &12395 &12892 &100 &126 &3716 &11702 &12055 \\
    $10^{-4\phantom{0}}$&131 &177 &4002 &14470 &14956 &83  &117 &2752 &9264  &9582  \\
    $10^{-6\phantom{0}}$&103 &135 &4386 &15035 &15409 &88  &132 &4170 &14421 &14774 \\
    $10^{-8\phantom{0}}$&73  &103 &3250 &11960 &12244 &101 &133 &3726 &13878 &14246 \\
    $10^{-10}$          &79  &109 &4088 &15527 &15825 &104 &139 &5378 &20291 &20674 \\
    \bottomrule
   \multicolumn{11}{c}{}
 \\\multicolumn{1}{c}{ }
  &\multicolumn{5}{c}{error-based termination}
  &\multicolumn{5}{c}{residual-based termination}
  \end{tabular}
\end{table}

The force term is $h(x_1, x_2) = - \sin(\omega x_1) \sin(\omega x_2)$,
with $\omega = \pi - \tfrac{1}{8}$. The control variable $z$
represents the Poisson diffusion coefficients that we are trying to
recover from the observed state $u_d$. We set $\alpha=10^{-4}$ as
the regularization parameter. The problem is almost identical to that
of \citet[\S 9.2]{EstrFrieOrbaSaun:2019a} but with an additional bound
constraint on the control variables (to ensure positivity of the
diffusion coefficients).

We discretize \eqref{eq:inverse-poisson} in two ways using $P_1$ finite elements on a uniform mesh of $1089$ (resp.\ $10201$) triangular elements and employ an identical discretization for the optimization variables
%\smarttodo{What is D? (Should be $\Omega$)}
$z \in L^2(\Omega)$, obtaining a problem with $n_z = 1089$ ($n_z = 10201$) controls and $n_u = 961$ ($n_u = 9801$) states, so that $n=n_u+n_z$. The control variables are discretized using piecewise linear elements. There are $m = n_u$ constraints, as we must solve the PDE on every interior grid point. For each problem, the target state is discretized on a finer mesh with 4 times more grid points and then interpolated onto the meshes previously described.

Although the problem on the smaller mesh was solved without the bound constraint in \cite[\S 9.2]{EstrFrieOrbaSaun:2019a}, the problem on the larger mesh could not be solved without explicitly enforcing the bound constraints because the control variables would go negative, causing the discretized PDE to be ill-defined.

We compute $x = (u,z)$ by applying KNITRO to~\eqref{eq:penalty-problem} with $\sigma = 10^{-2}$, using $B_2(x)$ as the Hessian approximation \eqref{eq:hess-approx-2} and initial point $u_0 = \ind$, $z_0 = \ind$. %We use KNITRO to optimize \eqref{eq:penalty-problem}. % with \LNLQ to (approximately) solve the symmetric form of \eqref{eq:aug-generic}.
We partition the Jacobian of the discretized constraints as $A(x)\T = \bmat{A_u(x)\T & A_z(x)\T}$, where
%\smarttodo{Dimensions are incompatible (Fixed)}
$A_u(x) \in \R^{n \times n}$, $A_z(x) \in \R^{m \times n}$ are the Jacobians for variables $u$, $z$ respectively. We use the preconditioner $\Pscr(x) = A_u(x)\T A_u(x)$, which amounts to performing two solves of a variable-coefficient Poisson equation (performed via direct solves). For this preconditioner, because the only bound constraints are $z \ge 0$, $Q(x) = \mbox{blkdiag}(I, Z)$ with $Z = \diag(z)$, so that
\begin{align*}
    \Pscr^{-1} A(x)\T Q(x) A(x) &= \Pscr^{-1} ( A_u(x)\T A_u(x) + A_z(x) Z A_z(x) )
\\  &= I + \Pscr^{-1} A_z(x) Z A_z(x).
\end{align*}
Thus $\sigma_{\min}(A(x) \Pscr^{-1/2}) \ge 1$, allowing us to bound the error via \LNLQ and to use both~\eqref{eq:terminate-error} and~\eqref{eq:terminate-residual} as termination criteria.

We choose $\epsilon = 10^{-8}$ in the stopping
conditions~\eqref{eq:experiment-stop}.
In
%\smarttodo{It's strange that (8.2b) is on the left and (8.2a) on the right (I flipped them)}
\cref{tab:inv-poisson} we vary
$\eta$, which defines the termination criteria of the linear system
solves~\eqref{eq:system-solve-accuracy}, and we record the number of
Hessian- and Jacobian-vector products. \cref{fig:inv_poisson} shows the target states and controls, and those that we recover on the two meshes (using \eqref{eq:terminate-error} and $\eta = 10^{-10}$).

We observed that for the smaller problem, KNITRO converged in a moderate number of outer iterations in all cases. With \eqref{eq:terminate-error}, we see that the number of Jacobian products tended to decrease as $\eta$ increased, except when $\eta=10^{-2}$ (for which the linear solves were too inaccurate). Using \eqref{eq:terminate-residual} showed a less clear trend. In cases with comparable outer iteration numbers, larger $\eta$ resulted in fewer Jacobian products. However, for moderate $\eta$ the number of outer iterations proved to be significantly smaller, resulting in a more efficient solve than when $\eta$ was too small or too large.

For the larger problem with termination condition
\eqref{eq:terminate-error}, the number of outer iterations
increased with $\eta$, the number of Lagrangian Hessian products
fluctuated somewhat, and Jacobian products tended to decrease. The
exception was $\eta = 10^{-8}$, which hit the sweet spot of solving
the linear systems sufficiently accurately to avoid many additional
outer iterations, but without performing too many iterations for each
linear solve. Using residual-based termination
\eqref{eq:terminate-residual} showed a less clear trend; Jacobian
products roughly decreased with increasing $\eta$ while the Hessian
products tended to oscillate. The sweet spot was hit with
$\eta=10^{-4}$, where the fewest outer iterations and operator
products were performed. For this problem, it appears that the
dependence of performance on the accuracy of the linear solves as
measured by the residual \eqref{eq:terminate-residual} is much more
nonlinear than when the linear solves are terminated according to the
error \eqref{eq:terminate-error}.

\subsection{2D Poisson-Boltzmann problem}

We now solve a control problem where the constraint is a 2D Poisson-Boltzmann equation:
%\smarttodo{Did you use two equation numbers intentionally? (No, my mistake)∫∫}
%\begin{alignat}{2}
\begin{equation}
\label{eq:poisson-boltzmann}
\begin{aligned}
    &\minimize{{u \in H_0^1(\Omega)},\ {z \in L^2(\Omega)}} &\quad& \half \int_\Omega \left( u - u_d \right)^2 \mathrm{d}x + \tfrac{1}{2} \alpha \int_\Omega z^2 \, \mathrm{d}x
\\  &\st &&
    \begin{array}[t]{rll}
         - \Delta u + \sinh(u) = & \hspace{-.5em} h + z & \mbox{in } \Omega,
      \\ u = & \hspace{-.5em} 0 & \mbox{on } \partial \Omega,
      \\z \ge & \hspace{-.5em} 0 & \mbox{in } \Omega.
    \end{array}
\end{aligned}
\end{equation}
%\end{alignat}
We use the same notation and $\Omega$ as in \cref{sec:inv-poisson}, with forcing term $h(x_1, x_2) = - \sin(\omega x_1) \sin(\omega x_2)$, $\omega = \pi - \tfrac{1}{8}$, and target state
\[
    u_d(x) = \begin{cases} 10 & \mbox{if } x\in [0.25,0.75]^2 \\
                           5  & \mbox{otherwise.}\end{cases}
\]
We discretized \eqref{eq:poisson-boltzmann} using $P_1$ finite elements on two uniform meshes with $1089$ (resp.\ $10201$) triangular elements, resulting in a problem with $n = 2050$ ($n = 20002$) variables and $m = 961$ ($m=9801$) constraints. The initial point was $u_0 = \ind$, $z_0 = \ind$.

\begin{figure}[t]
    \centering
    \begin{subfigure}{0.28\textwidth}
        \centering
        \includegraphics[height=103pt]{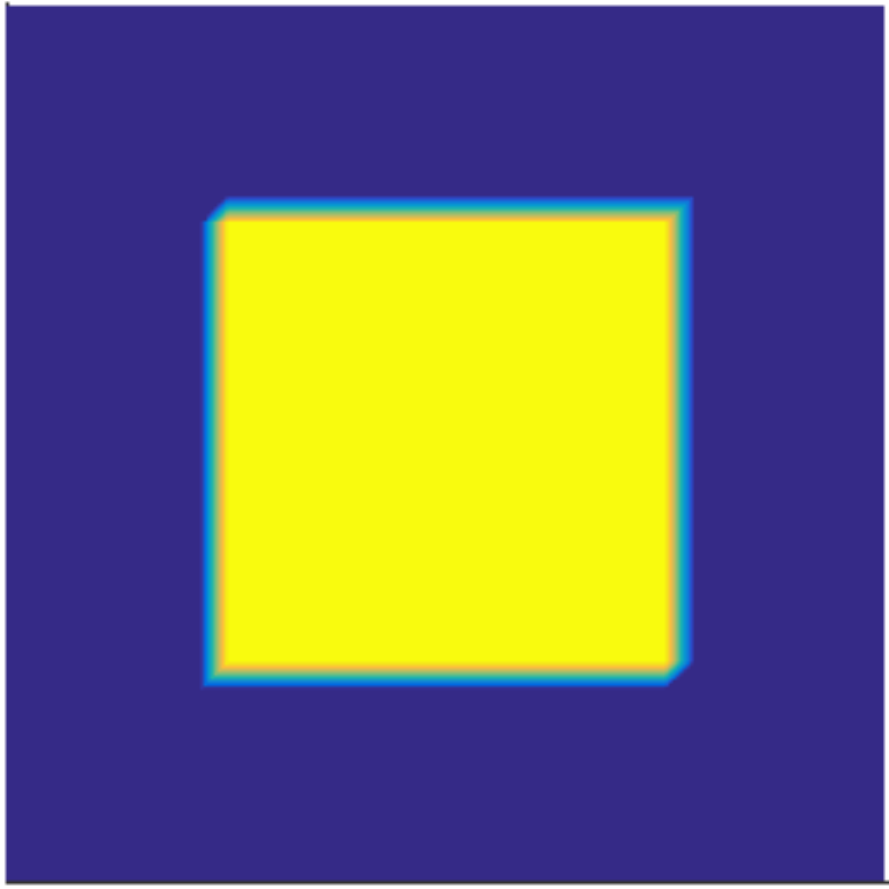}
        \caption{Target state $u_d$}
    \end{subfigure}
    ~
    \begin{subfigure}{0.28\textwidth}
        \centering
        \includegraphics[height=103pt,page=1]{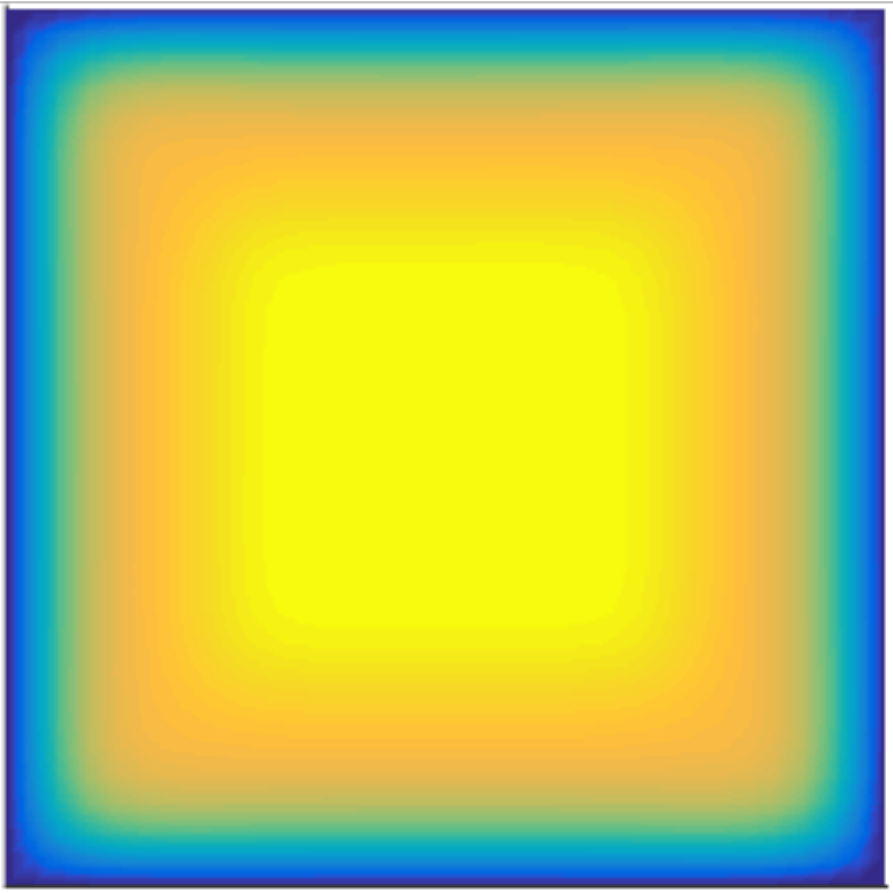}
        \caption{State (coarse mesh)}
    \end{subfigure}
    ~
    \begin{subfigure}{0.28\textwidth}
        \centering
        \includegraphics[height=103pt]{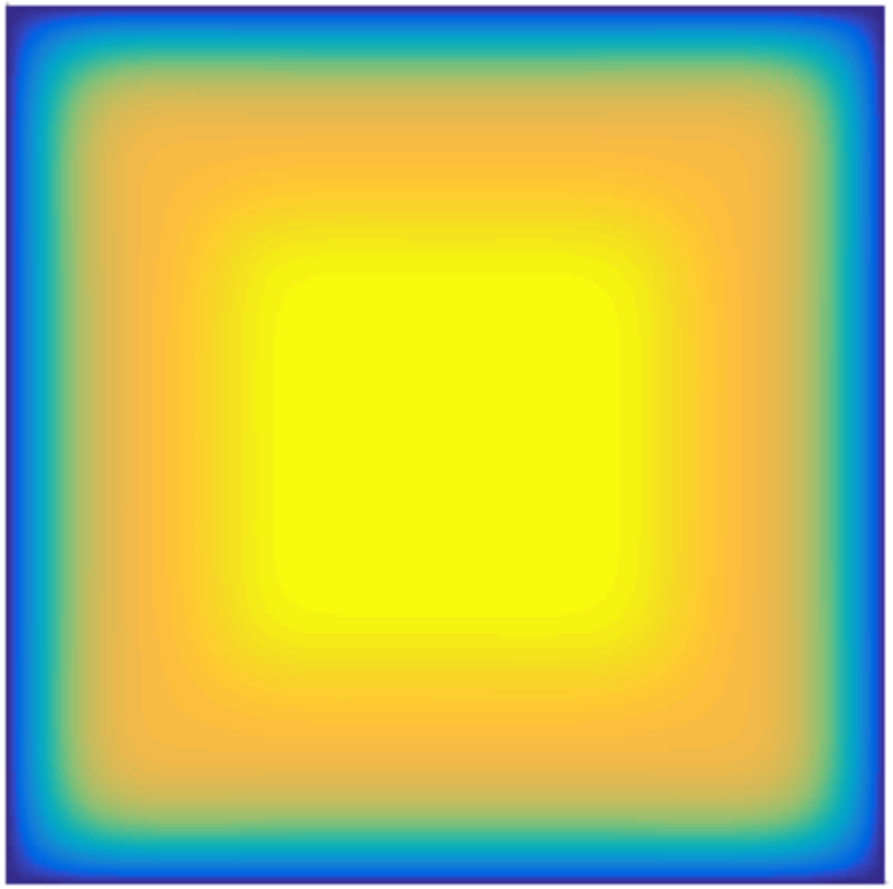}
        \caption{State (fine mesh)}
    \end{subfigure}
    ~
    \begin{subfigure}{0.08\textwidth}
        \centering
        \hspace{-11pt}
        \includegraphics[height=106pt,page=2]{ComputedStateSmall-2-down.pdf}
        \caption*{}
    \end{subfigure}
    \\
    \begin{subfigure}{0.28\textwidth}
        ~
    \end{subfigure}
    ~
    \begin{subfigure}{0.28\textwidth}
        \centering
        \includegraphics[height=103pt]{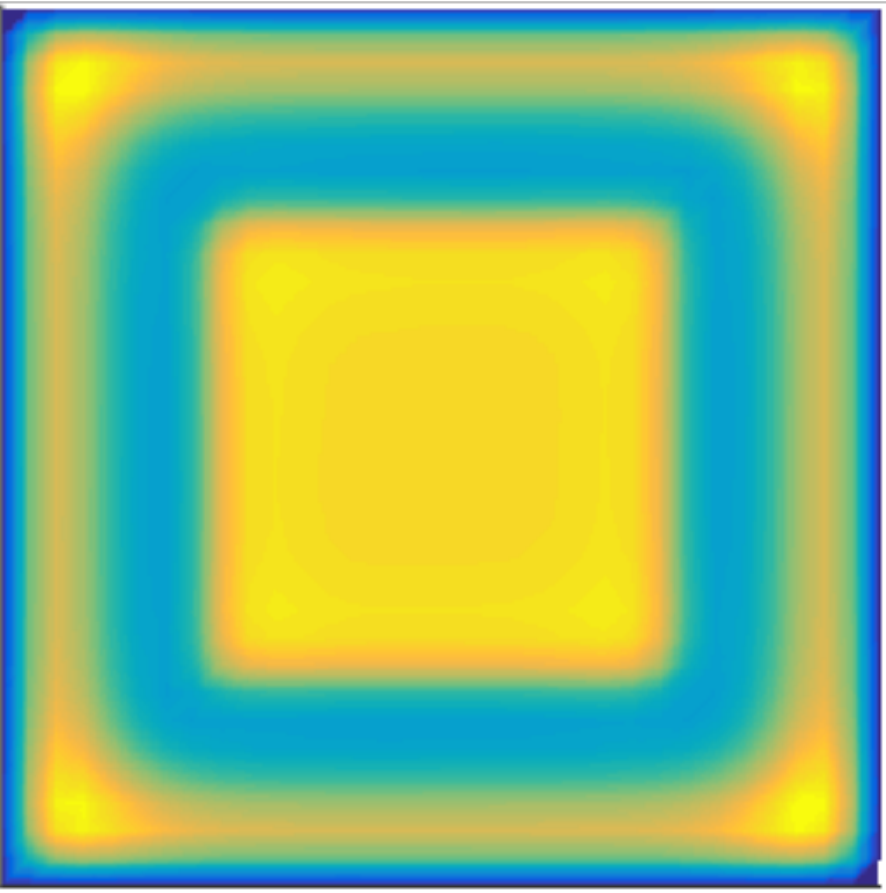}
        \caption{Control (coarse mesh)}
    \end{subfigure}
    ~
    \begin{subfigure}{0.28\textwidth}
        \centering
        \includegraphics[height=103pt]{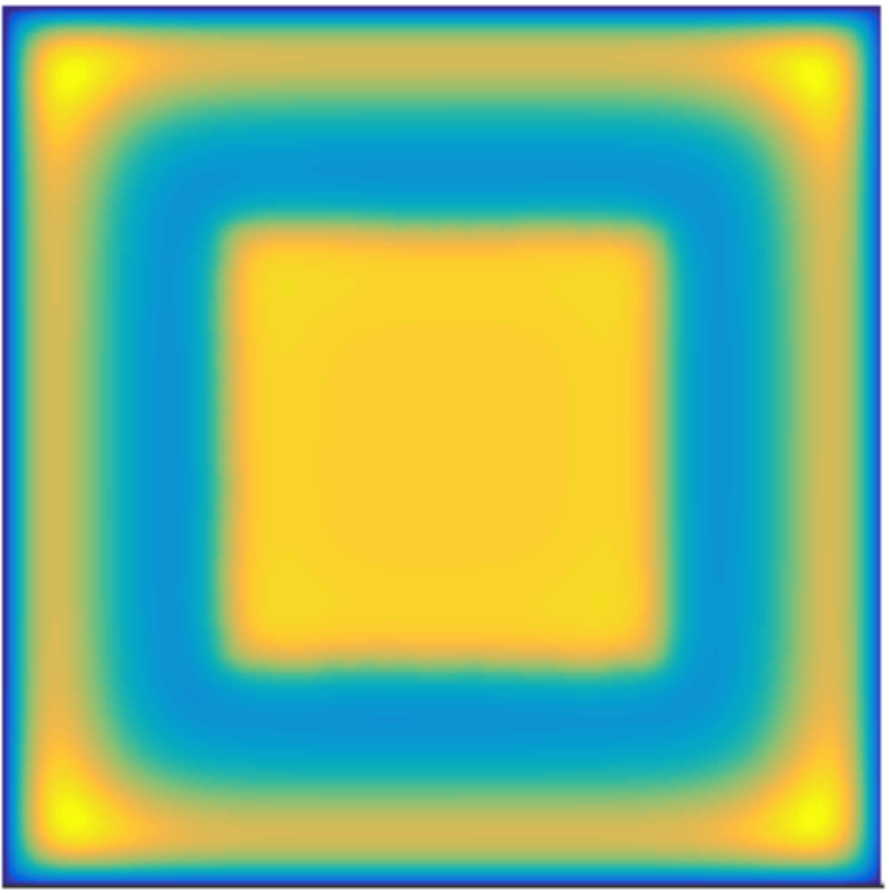}
        \caption{Control (fine mesh)}
    \end{subfigure}
    ~
    \begin{subfigure}{0.08\textwidth}
        \centering
        \includegraphics[height=103pt,page=2]{ComputedControlSmall-2-down.pdf}
        \caption*{}
    \end{subfigure}
    \caption{Target and computed states (top), and controls (bottom) for \eqref{eq:poisson-boltzmann}.}
    \label{fig:poisson-boltzmann}
\end{figure}

% Table 2
\begin{table}[H]
  \caption{Results from solving \eqref{eq:poisson-boltzmann} using KNITRO to optimize \eqref{eq:penalty-problem} with various $\eta$ in \eqref{eq:terminate-error} (left) and \eqref{eq:terminate-residual} (right) to terminate the %augmented
  linear system solves. The top (resp.\ bottom) table records results for the smaller problem with $n=2050$, $m=1089$ (resp.\ larger problem with $n=20002$, $m=10201$). We record the number of function/gradient evaluations ($\#f,g$), Lagrangian Hessian (\#$Hv$), Jacobian (\#($Av$), and adjoint Jacobian (\#$A\T v$) products.}
  \label{tab:poisson-boltzmann}
  \centering
  \begin{tabular}{@{} c |  c  c c  c  c | c c c c c @{}}
    \toprule
    $\eta$ & Its. & \#$f,g$ & \#$H v$ & \#$Av$ & \#$A\T v$ & Its. & \#$f,g$ & \#$H v$ & \#$Av$ & \#$A\T v$ \\
    \midrule
    $10^{-2\phantom{0}}$&19 &20 &1242 &3648 &3708 &19 &20 &1242 &3669 &3729 \\
    $10^{-4\phantom{0}}$&19 &20 &1252 &3753 &3813 &19 &20 &1244 &3762 &3822 \\
    $10^{-6\phantom{0}}$&19 &20 &1236 &3868 &3928 &19 &20 &1234 &3916 &3976 \\
    $10^{-8\phantom{0}}$&19 &20 &1244 &4169 &4229 &19 &20 &1236 &4286 &4346 \\
    $10^{-10}$          &19 &20 &1238 &4725 &4785 &19 &20 &1250 &4986 &5046 \\
    \midrule
    \multicolumn{11}{c}{} \vspace{-5pt}\\
    \toprule
    $10^{-2\phantom{0}}$&30 &37 &1524 &4426 &4531 &30 &37 &1524 &4468 &4573 \\
    $10^{-4\phantom{0}}$&30 &37 &1524 &4574 &4679 &30 &37 &1524 &4632 &4737 \\
    $10^{-6\phantom{0}}$&30 &37 &1524 &4813 &4918 &30 &37 &1558 &5033 &5138 \\
    $10^{-8\phantom{0}}$&30 &37 &1550 &5396 &5501 &30 &37 &1550 &5610 &5715 \\
    $10^{-10}$          &30 &37 &1550 &6224 &6329 &30 &37 &1558 &6582 &6687 \\
    \bottomrule
   \multicolumn{11}{c}{}
 \\\multicolumn{1}{c}{ }
  &\multicolumn{5}{c}{error-based termination}
  &\multicolumn{5}{c}{residual-based termination}
  \end{tabular}
\end{table}

We performed the same experiment as in \cref{sec:inv-poisson} using
%\smarttodo{How are you choosing the initial $\sigma$? (They were just chosen as the ones that work. How should we describe that?)}
$\sigma=10^{-1}$, and recorded the results in \cref{tab:poisson-boltzmann}. The target and computed state, and computed controls on the two meshes using \eqref{eq:terminate-error} with $\eta = 10^{-10}$ are given in \cref{fig:poisson-boltzmann}. We see that the results for both problems are more robust to changes in the accuracy of the linear solves. In all cases, the number of outer iterations and function/gradient evaluations were the same, and the number of Lagrangian Hessian products changed little. The number of Jacobian products steadily decreased with increasing $\eta$, with a 20--30\% drop in Jacobian products from $\eta=10^{-10}$ to $\eta=10^{-2}$.

% \pagebreak

\subsection{2D topology optimization}

\begin{figure}[t]
    \centering
    \begin{subfigure}{0.28\textwidth}
        \centering
        \includegraphics[height=95pt]{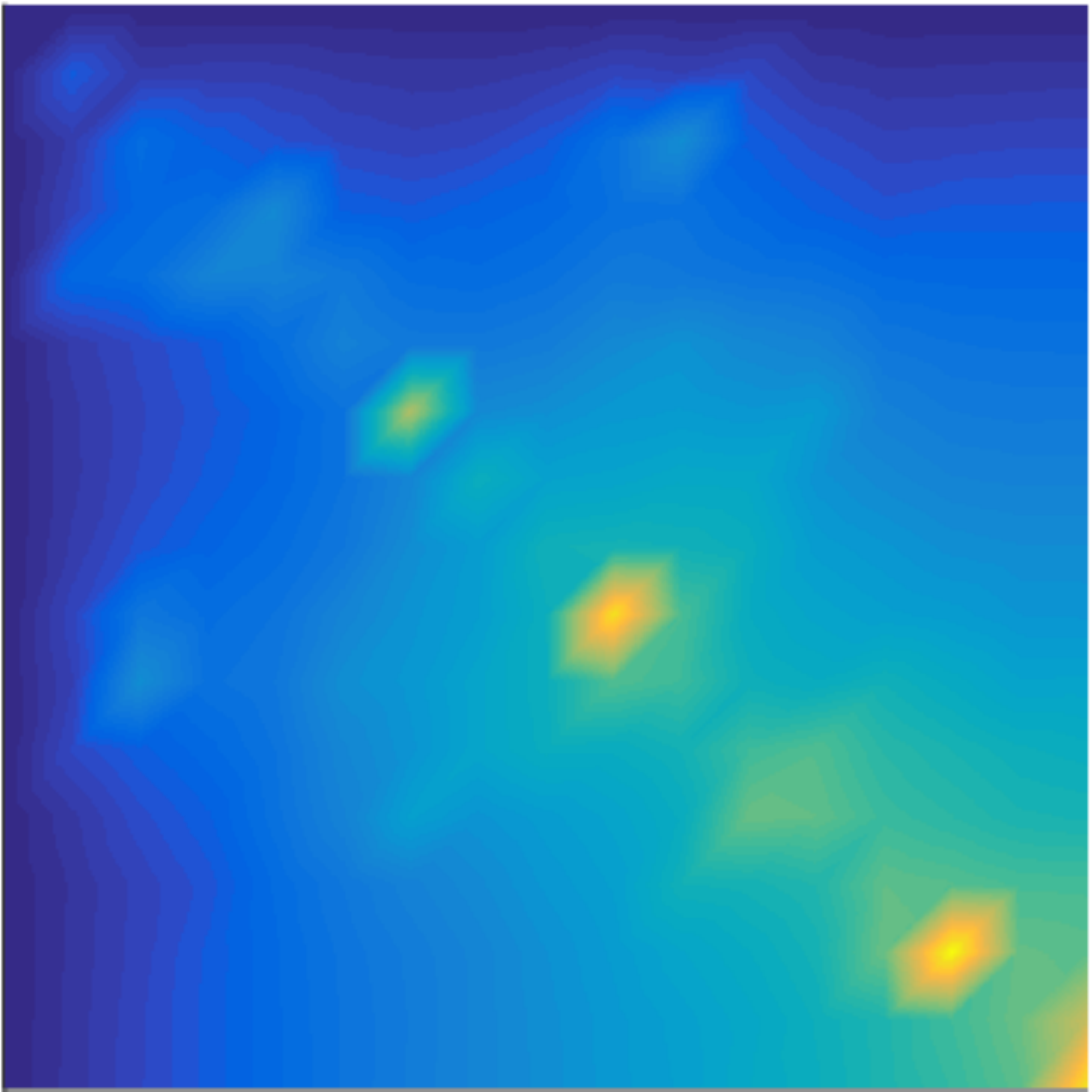}
        \caption{State (small mesh)}
    \end{subfigure}
    ~
    \begin{subfigure}{0.28\textwidth}
        \centering
        \includegraphics[height=95pt]{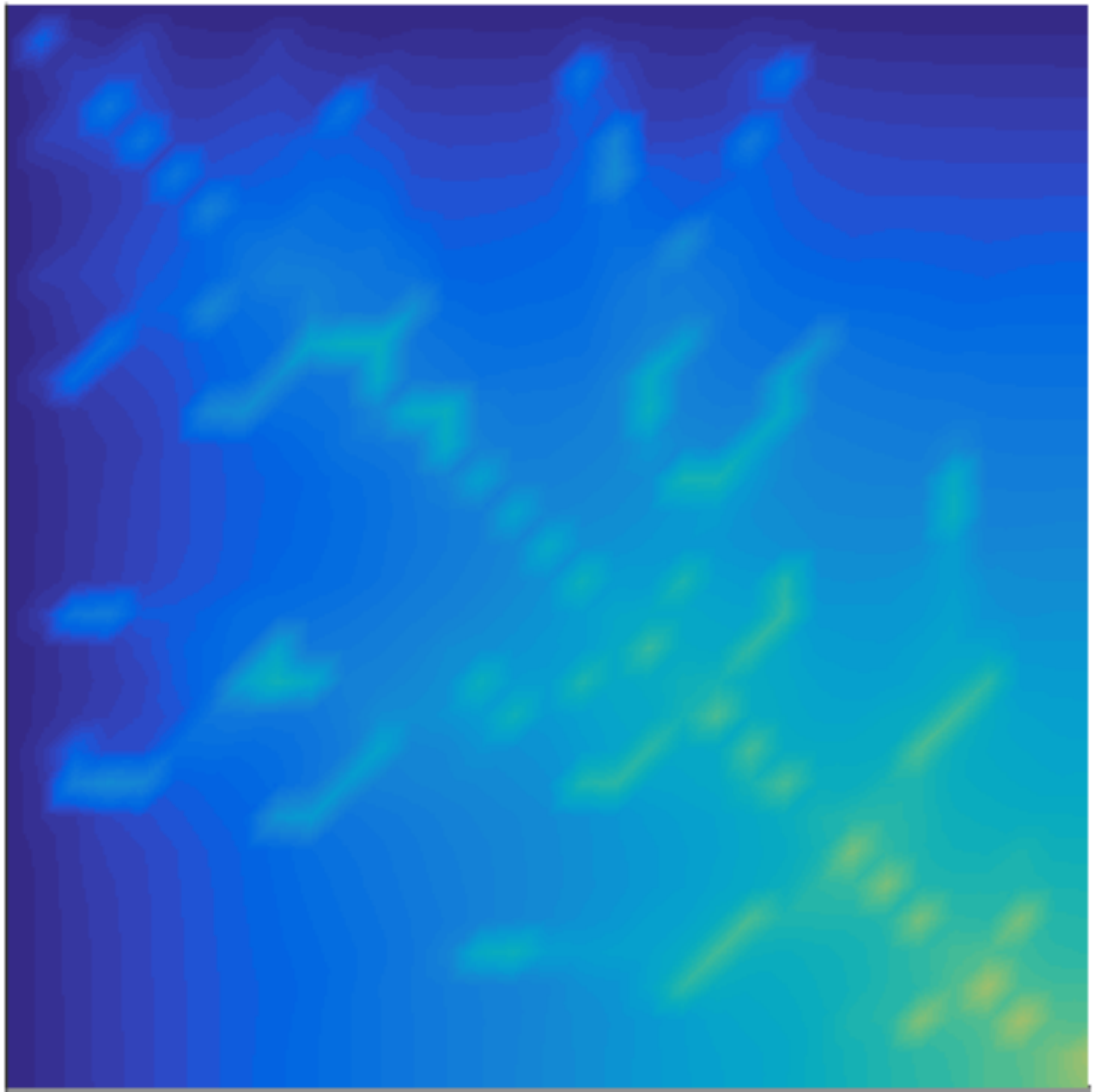}
        \caption{State (med mesh)}
    \end{subfigure}
    ~
    \begin{subfigure}{0.28\textwidth}
        \centering
        \includegraphics[height=95pt]{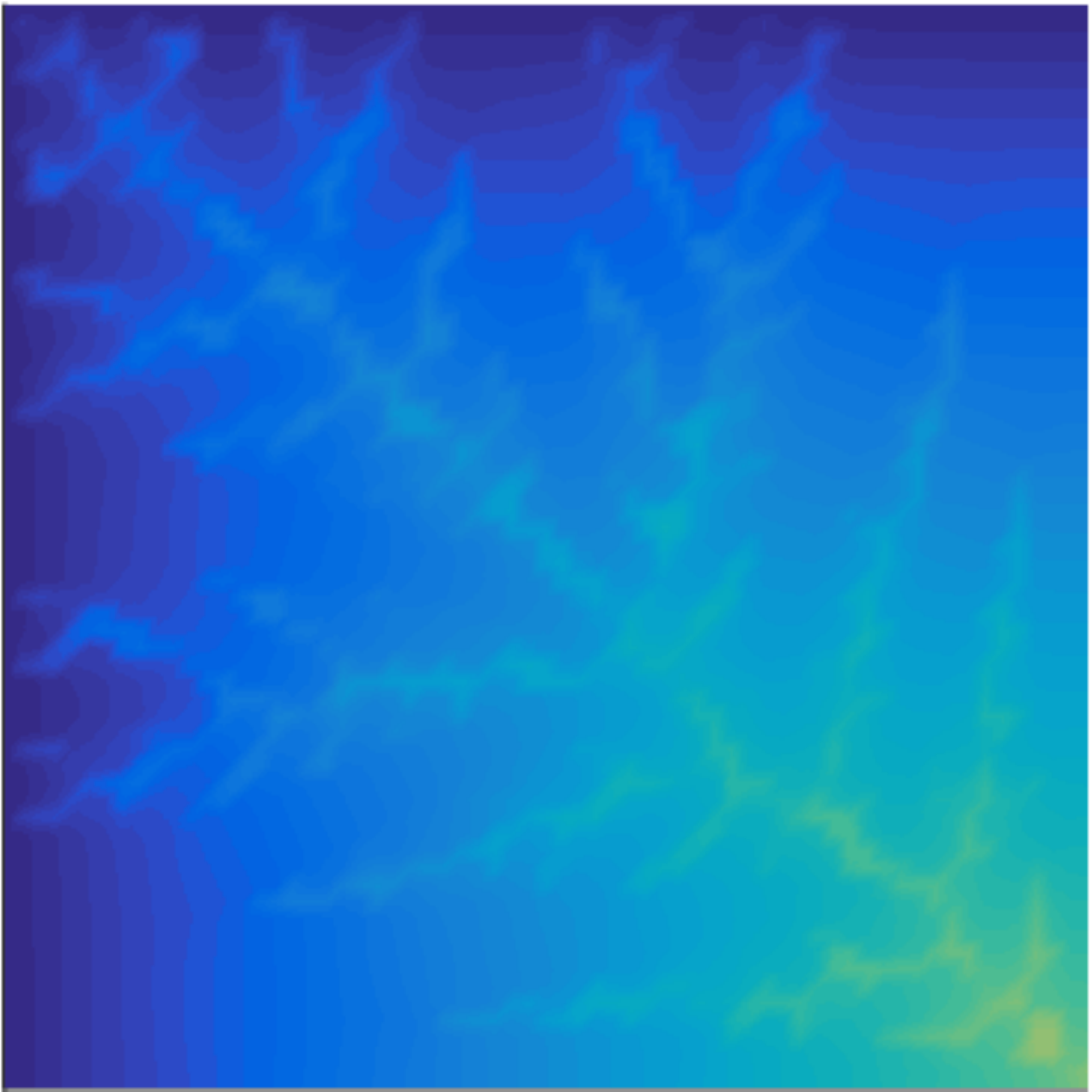}
        \caption{State (big mesh)}
    \end{subfigure}
    ~
    \begin{subfigure}{0.08\textwidth}
        \centering
        \vspace*{1pt}
        \includegraphics[height=96pt]{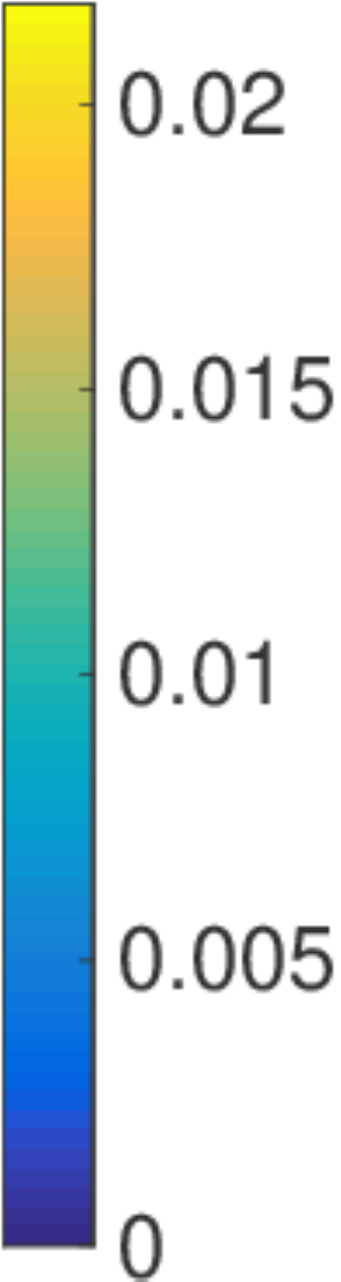}
        \caption*{}
    \end{subfigure}
    \\
    \begin{subfigure}{0.28\textwidth}
        \centering
        \includegraphics[height=95pt,page=1]{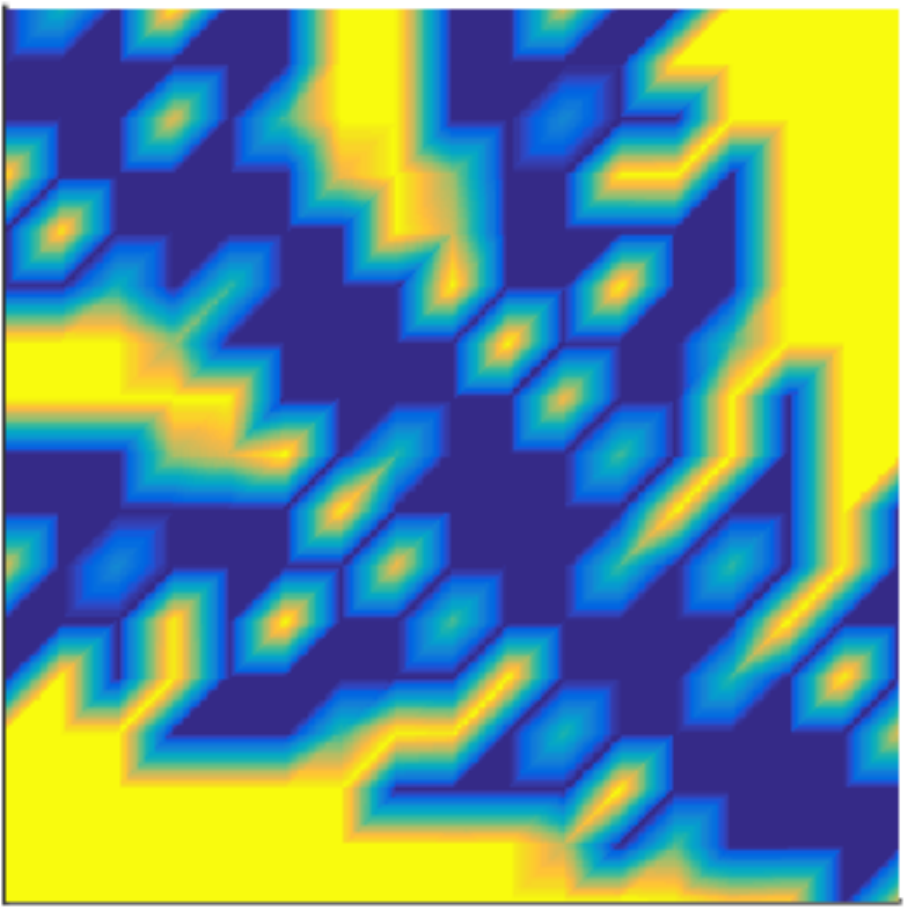}
        \caption{Control (small mesh)}
    \end{subfigure}
    ~
    \begin{subfigure}{0.28\textwidth}
        \centering
        \includegraphics[height=95pt,page=1]{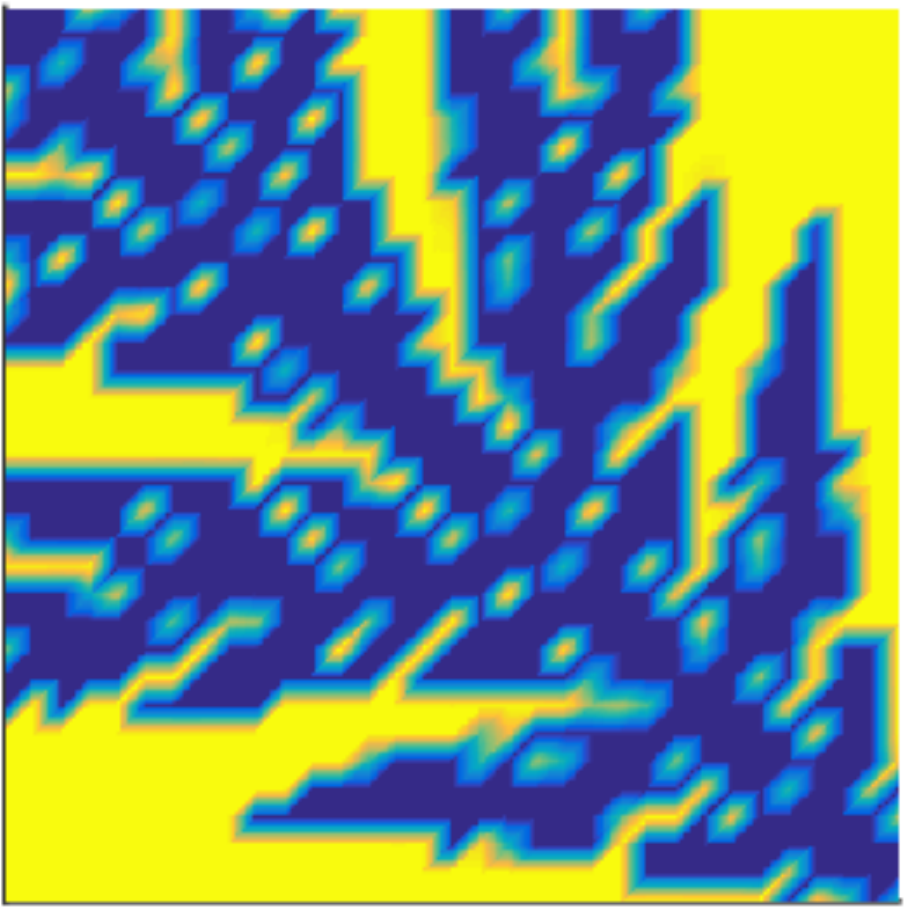}
        \caption{Control (med mesh)}
    \end{subfigure}
    ~
    \begin{subfigure}{0.28\textwidth}
        \centering
        \includegraphics[height=95pt,page=1]{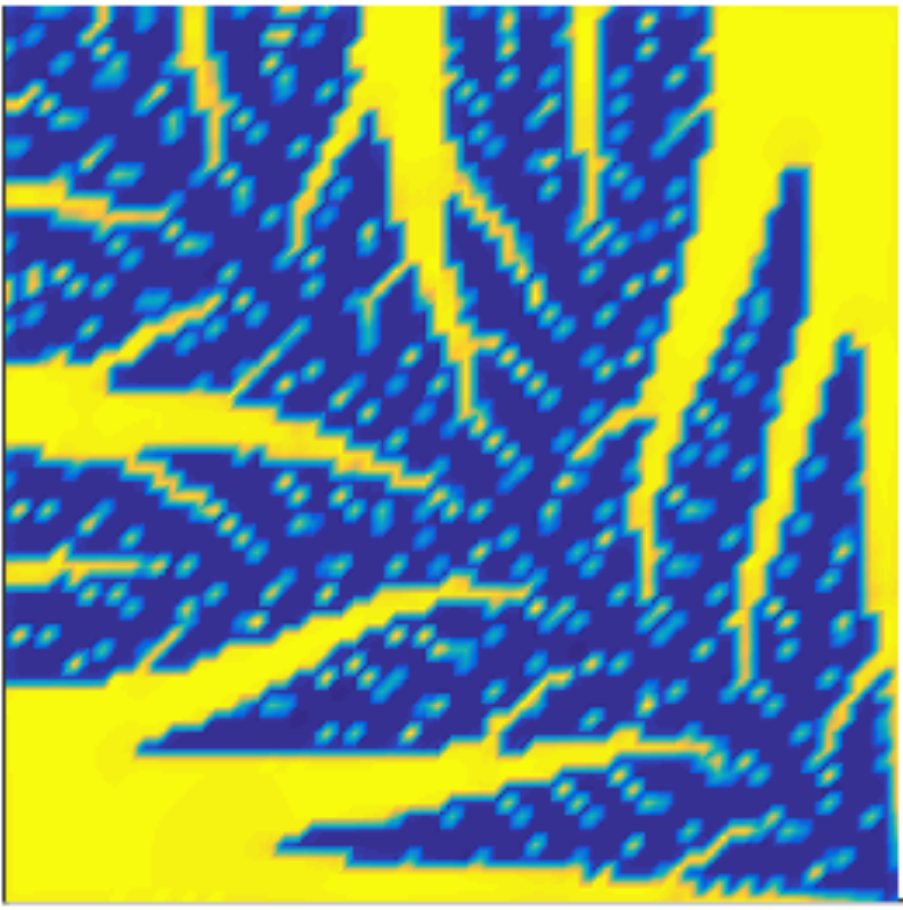}
        \caption{Control (big mesh)}
    \end{subfigure}
    ~
    \begin{subfigure}{0.08\textwidth}
        \centering
        \hspace*{-12pt}
        \includegraphics[height=96pt,page=2]{TopoptDensity16-2-down.pdf}
        \caption*{}
    \end{subfigure}
    \caption{Target and computed states (top), and controls (bottom) for \eqref{eq:poisson-boltzmann}.}
    \label{fig:topopt}
\end{figure}

We now solve the following 2D topology optimization problem from
\citet{GersBendSigm:2006}:

\begin{equation}
\label{eq:topology-problem}
\begin{aligned}
    \minimize{u \in H^1(\Omega), \, z \in L^2(\Omega)} \enspace & \int_\Omega f u \, \mathrm{d}x
\\  \st \quad &
    \begin{array}[t]{rll}
         \displaystyle{\textstyle{\int_{\Omega}} z \, \mathrm{d}x} \le & \hspace{-.5em} V &
      \\[2pt] - \nabla \cdot \left( k(z) \nabla u \right) = & \hspace{-.5em} f & \mbox{in } \Omega,
      \\ (k(z)\nabla u)\cdot\mathbf{n} = & \hspace{-0.5em} 0 &\mbox{on } \partial \Omega_1 = \{ (x,y) \mid x=0 \mbox{ or } y=1 \},
      \\ u = & \hspace{-.5em} 0 & \mbox{on } \partial \Omega_2 = \{ (x,y) \mid x=1 \mbox{ or } y=0 \},
      \\0 \le z \le & \hspace{-.5em} 1 & \mbox{in } \Omega,
    \end{array}
\end{aligned}
\end{equation}
where $k(z):\Omega \rightarrow \Omega$ defined by $k(z)(x) = 10^{-3} + (1-10^{-3}) z(x)^3$ for $x \in \Omega$, and $\mathbf{n}$ is the outward unit normal vector. The domain is $\Omega = [0,1]^2$, with load vector $f = 10^{-2}$, and $V = 0.4$. We discretize \eqref{eq:topology-problem} using finite elements as described by \cite{GersBendSigm:2006} on three grids: $16\times 16$, $32 \times 32$, and $64 \times 64$. This results in problems with $546$, $2114$, and $8321$ variables, and $256$, $1024$, and $4096$ equality constraints respectively. After discretization, we add a slack variable $s \ge 0$ for the first inequality constraint, so we have only equality constraints and bounds. The final problems then have one additional variable and constraint, with bounds on $z$ and $s$. %has $n = 8322$ variables and $m=4096$ constraints, with bounds on $z$ and $s$.

% Table 3
\begin{table}[t]
  \caption{Results from solving \eqref{eq:topology-problem} using KNITRO to optimize \eqref{eq:penalty-problem} with various $\eta$ in \eqref{eq:terminate-error} (left) and \eqref{eq:terminate-residual} (right) to terminate the %augmented
  linear system solves. Each table corresponds to a different mesh, with $16\times 16$ (top, $n=546$, $m=257$), $32\times 32$ (middle, $n=2114$, $m=1025$), and $64\times 64$ (bottom, $n=8322$, $m=4097$).
  We record the number of function/gradient evaluations ($\#f,g$), Lagrangian Hessian (\#$Hv$), Jacobian (\#($Av$), and adjoint Jacobian (\#$A\T v$) products. The symbol ``*" indicates that the problem failed to converge to a feasible point after 500 iterations.}
  \label{tab:topology}
  \centering
  \begin{tabular}{@{} c |  c  c  c  c  c | c c c c c @{}}
    \toprule
    $\eta$ & Its. & \#$f,g$ & \#$H v$ & \#$Av$ & \#$A\T v$ & Its. & \#$f,g$ & \#$H v$ & \#$Av$ & \#$A\T v$ \\
    \midrule
    $10^{-2\phantom{0}}$&176 &241 & 5296 & 15442 & 16342 &147 &230 & 3918 & 11462 & 12300 \\
    $10^{-4\phantom{0}}$&190 &286 & 6052 & 17694 & 18743 &171 &238 & 5634 & 16774 & 17660 \\
    $10^{-6\phantom{0}}$&164 &236 & 5266 & 15456 & 16329 &143 &199 & 3776 & 12019 & 12760 \\
    $10^{-8\phantom{0}}$&165 &239 & 5350 & 15743 & 16626 &176 &251 & 7100 & 23222 & 24152 \\
    $10^{-10}$          &185 &261 & 9096 & 26934 & 27903 &193 &289 & {\small 11420} & 39653 & 40714 \\
    \bottomrule
   \multicolumn{11}{c}{} \vspace{-5pt}\\
    \toprule
    $10^{-2\phantom{0}}$&219 &311 & 6598 & 19745 & 20898 &216 &319 & 6272 & 18381 & 19555 \\
    $10^{-4\phantom{0}}$&196 &265 & 5680 & 17073 & 18065 &189 &277 & 6382 & 18928 & 19949 \\
    $10^{-6\phantom{0}}$&190 &271 & 6190 & 15638 & 16642 &218 &302 & 7960 & 24383 & 25508 \\
    $10^{-8\phantom{0}}$&184 &272 & 4656 & 14050 & 15051 &211 &309 & 5868 & 19660 & 20799 \\
    $10^{-10}$          &184 &271 & 4396 & 13267 & 14265 &203 &291 & 5568 & 21526 & 22603 \\
    \bottomrule
   \multicolumn{11}{c}{} \vspace{-5pt}\\
    \toprule
    $10^{-2\phantom{0}}$&217 &340 & 4340 & 13966 & 15204 &* &* &* &* &* \\
    $10^{-4\phantom{0}}$&226 &348 & 4396 & 14068 & 15204 &* &* &* &* &* \\
    $10^{-6\phantom{0}}$&176 &272 & 3232 & 11218 & 12211 &191 &291 & 3508 & 18326 & 19391 \\
    $10^{-8\phantom{0}}$&185 &289 & 3356 & 11582 & 12635 &196 &296 & 3700 & 20888 & 21973 \\
    $10^{-10}$          &204 &298 & 4626 & 15412 & 16511 &190 &286 & 3480 & 23979 & 25028 \\
    \bottomrule
   \multicolumn{11}{c}{}
 \\\multicolumn{1}{c}{ }
  &\multicolumn{5}{c}{error-based termination}
  &\multicolumn{5}{c}{residual-based termination}
  \end{tabular}
\end{table}

We perform the same experiment as in \cref{sec:inv-poisson}, using $\sigma = 10^{-1}$ as the penalty parameter, and initial point $u_0 = \half V \ind$, $z_0 = \half V \ind$, $s_0 = V - \sum z_i = 0.2$. The linear constraint is kept explicit as in \cref{sec:explicit-constraints}. The results are recorded in \cref{tab:topology}. 

With \eqref{eq:terminate-error}, the number of outer iterations tends to increase with the mesh size; the trend is less clear with \eqref{eq:terminate-residual}. It is well known that such topology optimization problems become increasingly difficult numerically \citep{SigmPete:1998}, and typically require the use of a filter prior to solving the nonlinear optimization problem to improve its conditioning. Meshes refined as far as $128\times 128$ could not be solved directly using \eqref{eq:topology-problem}.

For a given mesh, when using \eqref{eq:terminate-error} the trend is like before: as $\eta$ increases the number of Jacobian products decreases (and in this case, so do the numbers of outer iterations and Lagrangian Hessian products), but this is only true until $\eta$ becomes too large and the linear solves become too coarse, causing slowed convergence. When \eqref{eq:terminate-residual} was used, we see a similar trend, except that when the linear solves are too coarse, KNITRO fails to converge.

\subsection{Explicit linear constraints}

We investigate the effect of maintaining the linear constraints
explicitly (\cref{sec:explicit-constraints}), using
%\smarttodo{There are just 2 problems; better say it}
some problems from the CUTEst test set \citep{GoulOrbaToin:2003} that have linear
constraints.
We use KNITRO to minimize $\phis$ with and without
linear constraints, because it can handle them explicitly. We use the corrected semi-normal equations to perform linear
solves, and
Hessian approximation
$B_1(x)$ \eqref{eq:hess-approx-1}. The
threshold penalty parameters \eqref{eq:14} and \eqref{eq:threshold-eq}
are computed from earlier optimal solutions when the linear
constraints were kept implicit
($\sigma^*_{\mbox{{\footnotesize impl}}}$) and explicit
($\sigma^*_{\mbox{{\footnotesize expl}}}$) respectively. The results
are recorded in \cref{tab:linear-constr}.

We observe that maintaining the linear constraints explicitly decreases the penalty parameter for all problems except \texttt{Channel400} ($\sigma^* = 0$ in both cases). KNITRO fails to find an optimal solution when the linear constraints are implicit and $\sigma < \sigma^*_{\mbox{{\footnotesize impl}}}$. This is because in the equality-constrained case $\phis$ is unbounded, and otherwise KNITRO stalls without converging to a feasible solution.
When $\sigma$ is sufficiently large, both versions converge (with and without explicit constraints); in most cases keeping the constraints requires fewer iterations, except for \texttt{Chain400}.
Although positive semidefiniteness of $\nabla^2 \phis(\xstar)$ is guaranteed in the relevant critical cone when $\sigma > \sigma^*$ (in either the implicit or explicit case), a larger value of $\sigma$ may sometimes be required because the curvature of $\phis$ away from the solution may be larger or ill-behaved.
%For example, for \texttt{prodpl0}, when $\sigma=3000$ is used, KNITRO converges in 472 iteration when linear constraints are implicit (for $\sigma=300$, KNITRO converges in 1067 iterations).

%We observe that maintaining the linear constraints explicitly decreases the penalty parameter for the \texttt{Chain} problems, and that KNITRO finds $\phis$ to be unbounded when $\sigma < \sigma^*_{\mbox{{\footnotesize impl}}}$ and the linear constraints are implicit. If $\sigma$ is large enough, both versions converge (with and without explicit linear constraints), and keeping the constraints implicit saves a few iterations.

For the \texttt{Channel} problems, the threshold parameter is zero in both cases. However, KNITRO converges quickly when the linear constraints are kept explicit, but otherwise fails to converge in a reasonable number of iterations. This phenomenon for the \texttt{Channel} problems appears to be independent of $\sigma$ (more values were investigated than are reported here). Even if the penalty parameter does not decrease, it appears beneficial to maintain some of the constraints explicitly.
%\smarttodo{Are the special cases that $\sigma^*$ is so small? Starting with such a small initial $\sigma$ would be unusual. Aren't there problems with more interesting thresholds that would also test our $\sigma$ update?}

\begin{table}[t]
  \caption{Results for problems with linear constraints
    (first three rows have only equality constraints).
    $m_{\mbox{{\scriptsize lin}}}$ and
    $m_{\mbox{{\scriptsize nln}}}$ are the number of
    linear and nonlinear constraints; $\sigma^*_{\mbox{{\scriptsize impl}}}$
    and $\sigma^*_{\mbox{{\scriptsize expl}}}$ are
    threshold penalty parameters when the linear constraints are
    handled implicitly and explicitly; $\sigma$ is the penalty
    parameter. The last two columns give the number of iterations
    before convergence; the symbol ``$*$'' indicates that unboundedness was detected,
    and ``-'' that 100 iterations were performed without converging.
    The solver exits when unboundedness is detected or an iterate
    satisfies~\eqref{eq:experiment-stop} with $\epsilon = 10^{-8}$.}
    % \label{tab:linear-constr}
    % \centering
    % \begin{tabular}{@{}c|c c c c c c c c@{}}
    %   \toprule
    %      Problem & $n$ & $m_{\mbox{{\footnotesize lin}}}$ & $m_{\mbox{{\footnotesize nln}}}$ & $\sigma^*_{\mbox{{\footnotesize impl}}}$ & $\sigma^*_{\mbox{{\footnotesize expl}}}$ & $\sigma$ & Impl. & Expl. \\
    %      \midrule
    %      \mr{\texttt{Chain100}} & \mr{202} & \mr{102} & \mr{1} & \mr{0.0047} & \mr{0} & $10^{-3}$ & $*$ & 13 \\
    %      &&&&&& 0.005 & 8 & 10 \\
    %      \hline
    %      \mr{\texttt{Chain200}} & \mr{402} & \mr{202} & \mr{1} & \mr{0.0024} & \mr{0} & $10^{-3}$ & $*$ & 11 \\
    %      &&&&&& 0.003 & 7 & 10 \\
    %      \hline
    %      \mr{\texttt{Chain400}} & \mr{802} & \mr{402} & \mr{1} & \mr{0.0012} & \mr{0} & $10^{-3}$ & $*$ & 10 \\
    %      &&&&&& 0.002  & 7 & 10 \\
    %      \hline
    %      \mr{\texttt{Channel100}} & \mr{800} & \mr{400} & \mr{400} & \mr{0} & \mr{0} & $10^{-3}$ & $-$ & 5 \\
    %      &&&&&& $1$ & $-$ & 5 \\
    %      \hline
    %      \mr{\texttt{Channel200}} & \mr{1600} & \mr{800} & \mr{800} & \mr{0} & \mr{0} & $10^{-3}$ & $-$ & 5 \\
    %      &&&&&& $1$ & $-$ & 5 \\
    %      \hline
    %      \mr{\texttt{Channel400}} & \mr{1600} & \mr{800} & \mr{800} & \mr{0} & \mr{0} & $10^{-3}$ & $-$ & 5 \\
    %              &&&&&& $1$ & $-$ & 5 \\
    %                                 \bottomrule
    % \end{tabular}
    \label{tab:linear-constr}
    \centering
    \small
    \begin{tabular}{@{}c|c c c c c c c c@{}}
      \toprule
         Problem & $n$ & $m_{\mbox{{\scriptsize lin}}}$ & $m_{\mbox{{\scriptsize nln}}}$ & $\sigma^*_{\mbox{{\scriptsize impl}}}$ & $\sigma^*_{\mbox{{\scriptsize expl}}}$ & $\sigma$ & Impl. & Expl. \\
         \midrule
         \mr{\texttt{Chain400}} & \mr{802} & \mr{402} & \mr{1} & \mr{0.0012} & \mr{0} & $10^{-3}$ & $*$ & 10 \\
         &&&&&& 0.002  & 7 & 10 \\
         \hline
         \mr{\texttt{Channel400}} & \mr{1600} & \mr{800} & \mr{800} & \mr{0} & \mr{0} & $10^{-3}$ & $-$ & 5 \\
         &&&&&& $1$ & $-$ & 5 \\
         \hline
         \mr{\texttt{hs113}} & \mr{18} & \mr{3} & \mr{5} & \mr{6.61} & \mr{3.39} & 6 & $*$ & 42 \\
         &&&&&& $7$ & 28 & 17 \\
         \hline
         \mr{\texttt{prodpl0}} & \mr{69} & \mr{25} & \mr{4} & \mr{211.9} & \mr{13.7} & $40$ & $-$ & 43 \\
         &&&&&& $300$ & $-$ %1067%
         & 30 \\
         \hline
         \mr{\texttt{prodpl1}} & \mr{69} & \mr{25} & \mr{4} & \mr{60.8} & \mr{3.56} & $10$ & $-$ & 22 \\
         &&&&&&  $70$ & 89 & 41 \\
         \hline
         \mr{\texttt{synthes3}} & \mr{38} & \mr{23} & \mr{19} & \mr{6.00} & \mr{0.66} & $2$ & $-$ & 12 \\
         &&&&&& $7$ & $35$ & 18 \\
         \bottomrule
    \end{tabular}
    \vspace*{-1\baselineskip}
\end{table}

\section{Discussion and concluding remarks}
\label{sec:discussion}

We derived a smooth extension of the penalty function by \cite{Fletcher:1970} as an extension to the implementation of \cite{EstrFrieOrbaSaun:2019a} to include bound constraints.
Our implementation is particularly promising for problems where augmented linear systems \eqref{eq:aug-generic} can be solved efficiently. We further demonstrated the merits of the approach on several PDE-constrained optimization problems.

Some limitations that are shared with the equality-constrained
case are avenues for future work. These include dealing with
the highly nonlinear nature of the penalty function, developing robust
penalty parameter updates and linear solve tolerance rules (for
inexact optimization solvers), preconditioning the trust-region
subproblems, and using cheaper second-derivative approximations (e.g.,
quasi-Newton updates) in conjunction with Hessian approximations
\eqref{eq:hess-approx-1}--\eqref{eq:hess-approx-2}. Possible
approaches for dealing with these issues are discussed by \citet[\S
10]{EstrFrieOrbaSaun:2019a}.

Bound constraints provide additional challenges for future work on top
of the equality-constrained case. For example, we would like to extend
the theory to problems with weaker constraint qualifications than
\ref{assump:licq}--\ref{assump:strict-complementarity}. A
regularization approach as in \citep[\S 6]{EstrFrieOrbaSaun:2019a} can
be employed when bound constraints are present, but it may need to be
refined to obtain similar convergence guarantees when
\ref{assump:licq} applies only at KKT points.

Another challenge is the possible numerical instability when iterates are close to the bounds, if the quantity $A(x)\T Q(x) A(x)$ becomes ill-conditioned. It would help to develop a specialized bound-constrained interior-point Newton-CG trust-region solver for \eqref{eq:penalty-problem} that carefully controls the distance to the bounds and attempts to minimize the number of approximate penalty Hessian products (as Hessian products are the most computationally intensive operation requiring two linear solves). We can also investigate other functions $Q(x)$ to approximate the complementarity conditions for KKT points, as different forms may have different advantages and limitations; for example, \eqref{eq:q-def} may cause premature termination if $\xstar$ is far from its bounds.

Our Matlab implementation can be found at \url{https://github.com/optimizers/FletcherPenalty}. To highlight the flexibility of Fletcher's approach, we implemented several options for applying various solvers to the penalty function and for solving the augmented systems, and other options discussed along the way.

\appendix
\section{Maintaining explicit constraints}
\label{app:explicit-constraints}

We discuss technical details about the penalty function when some of the constraints are linear and maintained explicitly as in~\eqref{eq:penalty-problem-exp}.
We define $\Ws(x) = \nabla \ws(x) \in \R^{n \times m_2}$, and $C(x) = \bmat{A(x) & B}$ as the Jacobian of all constraints. The operators $\gs(x)$, $\Hs(x)$, $S(x,v)$ and $T(x,w)$ are still defined over all constraints (e.g., $\gs(x) := g(x) - A(x)\ys(x) - B\ws(x)$), not just the nonlinear ones, and so they act on $C(x)$ and not just $A(x)$. Define
\begin{equation}
    \gsy(x) = g(x) - A(x) \ys(x)%, \qquad  \Hsy(x) = H(x) - \sum_{i=1}^m [\ys(x)]_i H_i(x)
\label{eq:partial-lag}
\end{equation}
as the gradient of the partial Lagrangian with respect to the nonlinear constraints $c(x)$ only (note that the linear constraints do not affect $\Hs$). The gradient and Hessian of the penalty function become
\begin{subequations} \label{eq:phi-grad-hess-eq}
  \begin{align}
    \nabla  \phis(x)  & = \gsy(x) - \Ys(x) c(x),
    \label{eq:phi-grad-eq}
    \\
    \nabla^2 \phis(x) & = \Hs(x) - A(x) \Ys(x)^T
                                    - \Ys(x) A(x)^T
                                    - \nabla_x \left[\Ys(x) c \right].
    \label{eq:phi-hess-eq}
  \end{align}
\end{subequations}

We restate the optimality conditions for \eqref{eq:nlp-exp} in terms of the penalty function. To do so, define the critical cones for \eqref{eq:nlp-exp} and \eqref{eq:penalty-problem-exp}, respectively, as
\[
    \widebar{\critcone}_{\phi}(\xstar, \zstar) = \critcone_{\phi}(\xstar,\zstar) \cap \{p \mid B\T p = 0\}, \qquad \widebar{\critcone}(\xstar, \zstar) = \critcone(\xstar,\zstar) \cap \{p \mid B\T p = 0\}.
\]

\begin{definition}[First-order KKT point]
\label{def:kkt-eq-1}
  A point $(\xstar, \zstar)$ is a
  first-order KKT point of~\eqref{eq:nlp-exp} if for any $\sigma \geq 0$ the
  following hold:
  \begingroup
  \allowdisplaybreaks
  \begin{subequations} \label{eq:1st-order-eq}
  \begin{align}
    \label{eq:100}
    \ell \le \xstar &\le u,
  \\c(\xstar)  &= 0, \label{eq:primal-feas-eq}
  \\B\T \xstar &= d,
  \\\nabla\phis(\xstar) &= B\wstar + \zstar,
  \\\zstar_j &= 0 \qquad \mbox{if } j \notin \Ascr(\xstar),
  \\\zstar_j &\ge 0 \qquad \mbox{if } \xstar_j = \ell_j, \label{eq:dual-ineq1-eq}
  \\\zstar_j &\le 0 \qquad \mbox{if } \xstar_j = u_j.  \label{eq:dual-ineq2-eq}
  \end{align}
  \end{subequations}
  \endgroup
  Then $\ystar:=\ys(\xstar)$ and $\wstar := \ws(\xstar)$ comprise the
  Lagrange multipliers of~\eqref{eq:nlp-exp} associated with $\xstar$. Note that by \ref{assump:strict-complementarity}, inequalities \eqref{eq:dual-ineq1-eq} and \eqref{eq:dual-ineq2-eq} are strict.
\end{definition}

\begin{definition}[Second-order KKT point]
\label{def:kkt-eq-2}
  The first-order KKT point
  $(\xstar,\zstar)$ satisfies the second-order necessary KKT condition
  for~\eqref{eq:nlp-exp} if for any $\sigma \geq 0$,
  \begin{equation} \label{eq:2nd-order-nec-eq}
  \begin{aligned}
    p^T \nabla^2\phis(\xstar) p \ge 0
    \quad
    \hbox{for all $p \in \widebar{\critcone}(\xstar, \zstar)$. }
  \end{aligned}
  \end{equation}
  The condition is sufficient if the inequality is strict.
\end{definition}

\begin{remark}
As before, if \eqref{eq:primal-feas-eq} is omitted, \cref{def:kkt-eq-1} defines first-order KKT points of \eqref{eq:penalty-problem-exp}. Similarly, replacing $\widebar{\critcone}(\xstar,\zstar)$ by $\widebar{\critcone}_{\phi}(\xstar,\zstar)$ in \cref{def:kkt-eq-2} defines second-order KKT points of \eqref{eq:penalty-problem-exp}.
\end{remark}

\subsection{Proof of \texorpdfstring{\cref{thm:threshold-eq}}{Theorem 6}}
\label{sec:proof-thm-threshold}

Observe that the multiplier estimates $\ys(x)$ and $\ws(x)$ satisfy
\begin{equation}
    \label{eq:37}
    C(x)^T Q(x) C(x) \bmat{\ys(x) \\ \ws(x)} = C(x)^T Q(x) g(x) - \sigma \bmat{c(x) \\ B\T x - d}.
\end{equation}

Proof of \eqref{eq:34}: We drop the argument $x$ from operators and assume that all are evaluated at~$\xbar$. Because $\xbar$ is a first-order KKT point for \eqref{eq:penalty-problem-exp}, we need only show that $c(\xbar) = 0$. Further,  $Q (\nabla \phis - B \wstar) = 0$ at $\xbar$, or equivalently,
\begin{align*}
    Q B\wstar = Q \left( g - A \ys - \Ys c \right).% Q(\xbar) B\wstar = Q(\xbar) \left( g(\xbar) - A(\xbar) \ys(\xbar) - \Ys(\xbar) c(\xbar) \right).
\end{align*}
Multiplying both sides by $C\T$ and using \eqref{eq:37} we have
% \begin{align*}
%     \bmat{A(\xbar)\T Q(\xbar) B\wstar \\ B\T Q(\xbar) B \wstar} = \sigma \bmat{c(\xbar) \\ 0} + \bmat{A(\xbar)\T Q(\xbar) B \ws(\xbar) \\ B\T Q(\xbar) B \ws(\xbar)} - \bmat{A(\xbar)\T Q(\xbar) \Ys(\xbar) c(\xbar) \\ B\T Q(\xbar) \Ys(\xbar) c(\xbar)},
% \end{align*}
\begin{align*}
    \bmat{A\T Q B\wstar \\ B\T Q B \wstar} = \sigma \bmat{c \\ 0} + \bmat{A\T Q B \ws \\ B\T Q B \ws} - \bmat{A\T Q \Ys c \\ B\T Q \Ys c},
\end{align*}
so that $\ws = \wstar + (B\T Q B)^{-1} B\T Q \Ys c$. %$\ws(\xbar) = \wstar + (B\T Q(\xbar) B)^{-1} B\T Q(\xbar) \Ys(\xbar) c(\xbar)$.
Substituting $\ws(\xbar)$ into the first block of equations and rearranging gives
% \begin{equation*}
%      A(\xbar) Q(\xbar)^{1/2} \Pbar_{Q(\xbar)^{1/2}B} Q(\xbar)^{1/2} \Ys(\xbar) c(\xbar) = \sigma c(\xbar).
% \end{equation*}
\begin{equation*}
     A Q^{1/2} \Pbar_{Q^{1/2}B} Q^{1/2} \Ys c = \sigma c.
\end{equation*}
The triangle inequality gives $\sigma \norm{c} \le \norm{A^T Q^{1/2} \Pbar_{Q^{1/2}B} Q^{1/2} \Ys} \norm{c}$, implying $c = 0$. Then $\ws = \wstar$ and $\xbar$ is a first-order KKT point for \eqref{eq:nlp-exp}.

Proof of \eqref{eq:36}: As in the proof of \eqref{eq:16}, we differentiate \eqref{eq:37} to obtain
% \begin{equation}
%     \label{eq:38}
%     C(x)\T C(x) \bmat{\Ys(x)^T \\ \Ws(x)^T} = C(x)^T [\Hs(x) - \sigma I] + S(x, \gs(x)).
% \end{equation}
\begin{equation}
    \label{eq:38}
    \begin{aligned}
    &C(x)\T Q(x) C(x) \bmat{\Ys(x)^T \\ \Ws(x)^T}
    \\ &\qquad = C(x)\T \left[ Q(x)\Hs(x) - \sigma I + R(x, \gs(x)) \right] + S(x,Q(x)\gs(x)).
    \end{aligned}
\end{equation}
For the remainder of the proof, we assume all operators are evaluated at $\xstar$.
Because $\xstar$ satisfies first-order conditions \eqref{eq:1st-order-eq}, $Q \gs = 0$ independently of $\sigma$, so $S(Q\gs) = 0$. Let $P_{Q^{1/2}C} := P_{Q^{1/2} C(\xstar)}(\xstar)$, so that from \eqref{eq:38} we have
\begin{equation}
    \label{eq:39}
    Q^{1/2} \left( A \Ys^T + B \Ws^T \right) Q^{1/2} = P_{Q^{1/2}C} \left[Q^{1/2} \Hs Q^{1/2} - \sigma I + R(\gs) Q^{1/2}\right].
\end{equation}
Observe that if $p\in \widebar{\critcone}_{\phi}(\xstar,\zstar)$, then $p = Q^{1/2} \pbar$ for some $\pbar \in \widebar{\critcone}_{\phi}(\xstar,\zstar)$. Because $Q^{1/2} \gs = 0$, we have $R(\gs) p = 0$.

Substituting \eqref{eq:39} into \eqref{eq:phi-hess-eq}, % using $\Hs(\xstar) = \hLag(\xstar,\ystar)$,
and $P_{Q^{1/2}C} + \Pbar_{Q^{1/2}C} = I$ gives
\begin{align*}
    &p^T \nabla^2 \phis(\xstar) p \ge 0
\\  \iff& \pbar^T Q^{1/2} \left( \Hs - A \Ys\T - \Ys A^T\right) Q^{1/2} \pbar \ge 0
\\  \iff& \pbar^T \left( \Pbar_{Q^{1/2}C} Q^{1/2} \Hs Q^{1/2} \Pbar_{Q^{1/2}C} - P_{Q^{1/2}C} Q^{1/2} \Hs Q^{1/2} P_{Q^{1/2}C} + 2\sigma P_{Q^{1/2}C}  \right) \pbar 
\\ &\qquad - p\T \left( B\Ws\T + \Ws B\T  \right) p \ge 0.
\end{align*}
% \begin{equation*}
%     \nabla^2 \phis(\xstar_2) = \Pbar_C \hLag(\xstar_2,\ystar) \Pbar_C - P_C \hLag(\xstar_2,\ystar) P_C + 2 \sigma P_C - B \Ws(x)^T - \Ws(x) B^T .
% \end{equation*}
Because $\Hs(\xstar) = \hLag(\xstar,\ystar)$, $0 = B^T p = B^T Q^{1/2} \pbar$, we can write $\pbar = \Pbar_{\Bbar} q$ with $\Bbar = Q^{1/2}B$ and hence
\begin{align*}
    0 \leq& p^T \nabla^2 \phis(\xstar_2) p
\\  \iff 0 \preceq& \Pbar_{\Bbar} \Pbar_{Q^{1/2}C} \hLag(\xstar,\ystar) \Pbar_{Q^{1/2}C}\Pbar_{\Bbar}  
\\ &\qquad - \Pbar_{\Bbar} P_{Q^{1/2}C} \hLag(\xstar,\ystar) P_{Q^{1/2}C} \Pbar_{\Bbar} + 2      \sigma \Pbar_{\Bbar} P_{Q^{1/2}C} \Pbar_{\Bbar}.
\end{align*}
As before, the first term is positive semi-definite, so we only need that
\[
- \Pbar_{\Bbar} P_{Q^{1/2}C} \hLag(\xstar,\ystar) P_{Q^{1/2}C} \Pbar_{\Bbar} + 2 \sigma \Pbar_{\Bbar} P_{Q^{1/2}C} \Pbar_{\Bbar} \succeq 0,
\]
which is equivalent to $\sigma \geq \sigmabar$. \qed

\subsubsection{Evaluating the penalty function and derivatives}
\label{sec:eval-explicit}

We again drop the arguments on functions and assume they are evaluated at a point $x$ for some $\sigma$:
\[
y = \ys(x), \quad A = A(x), \quad \Ys = \Ys(x), \quad
H_\sigma = H_\sigma(x), \quad S_\sigma = S_\sigma(x,g_\sigma(x)),
\quad \hbox{etc.}
\]
We focus on the nonsymmetric linear systems; the corresponding symmetric linear systems can be derived similarly to \cref{sec:comp-penalty-funct}.

The multipliers for evaluating the penalty function are obtained by solving
\begin{align}
    \label{eq:aug-mult-eq}
    \bmat{I & A & B \\ A\T Q &  & \\ B\T Q & & } \bmat{\gs \\ \ys \\ \ws} = \bmat{ g \\ \sigma c \\ \sigma (Bx - d)}.
\end{align}

To compute the gradient and Hessian products, we use the identity
\begin{equation}
    C\T Q C \bmat{\Ys^T \\ \Ws^T}  = C\T \left[ Q\Hs - \sigma I + \Rs \right] + \Ss
\end{equation}
to obtain the necessary products with $\Ys$ and $\Ys\T$. Observe that
\begin{align*}
  \Ys u = \bmat{\Ys & \Ws}\bmat{u \\ 0}, \qquad \Ys\T v = \bmat{I & 0} \bmat{\Ys^T \\ \Ws^T} v,
\end{align*}
so that \cref{alg:Yu-sym} and \cref{alg:YTv-sym} can be applied.

Note that to compute the gradient in~\eqref{eq:phi-grad-eq}, $\gsy$ is not available directly from the solution to \eqref{eq:aug-mult-eq} and must be computed explicitly using~\eqref{eq:partial-lag}.

Approximate products with $\nabla^2 \phis$ can be computed via
\begin{align*}
    \nabla^2\phis &\approx B_1 := \Hs - A\Ys^T - \Ys A^T %\label{eq:hess-approx-1-eq}
\\  &\phantom{\approx B_1 :}= \Hs - \bmat{A & 0} (C\T Q C)^{-1} C\T (Q \Hs - \sigma I  + \Rs) - \bmat{A & 0} (C\T Q C)^{-1} S_{\sigma} %\nonumber
\\ &\phantom{\approx B_1 :=} \qquad - (\Hs Q - \sigma I + \Rs) C (C\T Q C)^{-1} \bmat{A^T \\ 0} - S_{\sigma} (C\T Q C)^{-1} \bmat{A^T \\ 0} %\nonumber
\\  &\approx B_2 := \Hs - \bmat{A & 0} (C\T Q C)^{-1} C\T (Q \Hs - \sigma I  + \Rs) %\nonumber
\\ &\phantom{\approx B_1 :=} \qquad - (\Hs Q - \sigma I + \Rs) C (C\T Q C)^{-1} \bmat{A^T \\ 0}. %\nonumber
\end{align*}
For products with the weighted-pseudoinverse and its transpose, we can compute
\[
\bmat{u_1 \\ u_2} = (C\T Q C)^{-1} C\T v, \qquad v = C (C\T Q C)^{-1} \bmat{u_1 \\ u_2}
\]
by solving the respective block systems
\begin{equation}
    \label{eq:aug-pseudo}
    \bmat{I & QA & QB \\ A\T &  & \\ B\T & & } \bmat{t \\ u_1 \\ u_2} = \bmat{v \\ 0 \\ 0}, \qquad \bmat{I & A & B \\ A\T Q &  & \\ B\T Q & & } \bmat{v \\ t_1 \\ t_2} = \bmat{0 \\ -u_1 \\ -u_2}.
\end{equation}
Thus we can obtain the same types of Hessian approximations as \eqref{eq:hess-approx}, again with two augmented system solves per product.

\section*{Acknowledgements}
We would like to express our deep gratitude to Drew Kouri for supplying PDE-constrained optimization problems in Matlab, for helpful discussions throughout this project, and for hosting the first author for two weeks at Sandia National Laboratories.
We are also grateful to the reviewers for their careful reading and many helpful questions and suggestions.

\small
% the only difference between myabbrvnat and the standard abbrvnat is that
% the DOI field for technical report is printed
\bibliographystyle{myabbrvnat}  %{siamplain} %
\bibliography{paper}

\providecommand{\noopsort}[1]{}
\begin{thebibliography}{32}
\providecommand{\natexlab}[1]{#1}
\providecommand{\url}[1]{\texttt{#1}}
\expandafter\ifx\csname urlstyle\endcsname\relax
  \providecommand{\doi}[1]{doi: #1}\else
  \providecommand{\doi}{doi: \begingroup \urlstyle{rm}\Url}\fi

\bibitem[Anitescu(2000)]{Anit:2000}
M.~Anitescu.
\newblock On solving mathematical programs with complementarity constraints as
  nonlinear programs.
\newblock Preprint ANL/MCS-P864-1200, Argonne National Laboratory, 2000.

\bibitem[Arioli(2013)]{Ario:2013}
M.~Arioli.
\newblock Generalized {G}olub-{K}ahan bidiagonalization and stopping criteria.
\newblock \emph{SIAM J. Matrix Anal. Appl.}, 34\penalty0 (2):\penalty0
  571--592, 2013.
\newblock \doi{10.1137/120866543}.

\bibitem[Bertsekas(1975)]{Bert:1975}
D.~P. Bertsekas.
\newblock Necessary and sufficient conditions for a penalty method to be exact.
\newblock \emph{Math. Program.}, 9:\penalty0 87--99, 1975.

\bibitem[Bertsekas(1982)]{Bert:1982}
D.~P. Bertsekas.
\newblock \emph{Constrained Optimization and {L}agrange Multiplier Methods}.
\newblock Academic Press, New York, 1982.

\bibitem[Bj\"orck and Paige(1994)]{PaigBjor:1994}
A.~Bj\"orck and C.~C. Paige.
\newblock Solution of augmented linear systems using orthogonal factorizations.
\newblock \emph{BIT}, 34\penalty0 (1):\penalty0 1--24, 1994.
\newblock \doi{10.1007/BF01935013}.

\bibitem[Boggs et~al.(1992)Boggs, Tolle, and Kearsley]{BoggTollKear:1991}
P.~T. Boggs, J.~W. Tolle, and A.~J. Kearsley.
\newblock A merit function for inequality constrained nonlinear programming
  problems.
\newblock Internal Report {NISTIR} 4702, {A}pplied and {C}omputational
  {M}athematics {D}ivision, {N}ational {I}nstitute of {S}tandards and
  {T}echnology, {G}aithersburg, {MD}, {USA}, 1992.

\bibitem[Byrd et~al.(2006)Byrd, Nocedal, and Waltz]{ByrdNoceWalt:06}
R.~H. Byrd, J.~Nocedal, and R.~A. Waltz.
\newblock {KNITRO}: An integrated package for nonlinear optimization.
\newblock In G.~di~Pillo and M.~Roma, editors, \emph{Large-Scale Nonlinear
  Optimization}, pages 35--59. Springer-Verlag, New York, 2006.

\bibitem[Chen(2000)]{Chen:2000}
X.~Chen.
\newblock Smoothing methods for complementarity problems and their
  applications: a survey.
\newblock \emph{J. Oper. Res. Soc. Japan}, 43\penalty0 (1):\penalty0 32--47,
  2000.
\newblock ISSN 0453-4514.
\newblock \doi{10.1016/S0453-4514(00)88750-5}.
\newblock New trends in mathematical programming (Kyoto, 1998).

\bibitem[Conn et~al.(2000)Conn, Gould, and Toint]{ConnGoulToin:2000}
A.~R. Conn, N.~I.~M. Gould, and {\mbox{Ph}}.~L. Toint.
\newblock \emph{Trust-Region Methods}.
\newblock {MPS-SIAM} Series on Optimization. SIAM, Philadelphia, 2000.

\bibitem[Craig(1955)]{Craig:1955}
J.~E. Craig.
\newblock The {N}-step iteration procedures.
\newblock \emph{Journal of Mathematics and Physics}, 34\penalty0 (1):\penalty0
  64--73, 1955.

\bibitem[Di~Pillo and Grippo(1984)]{PilloGrippo:1984}
G.~Di~Pillo and L.~Grippo.
\newblock A class of continuously differentiable exact penalty function
  algorithms for nonlinear programming problems.
\newblock In E.~P. Toft-Christensen, editor, \emph{System Modelling and
  Optimization}, page 246–256. Springer-Verlag, Berlin, 1984.

\bibitem[Di~Pillo and Grippo(1985)]{PillGrip:1985}
G.~Di~Pillo and L.~Grippo.
\newblock A continuously differentiable exact penalty function for nonlinear
  programming problems with inequality constraints.
\newblock \emph{SIAM J. Control Optim.}, 23\penalty0 (1):\penalty0 72--84,
  1985.
\newblock ISSN 0363-0129.
\newblock \doi{10.1137/0323007}.

\bibitem[Estrin and Greif(2018)]{EstrGrei:2018}
R.~Estrin and C.~Greif.
\newblock S{PMR}: a family of saddle-point minimum residual solvers.
\newblock \emph{SIAM J. Sci. Comput.}, 40\penalty0 (3):\penalty0 A1884--A1914,
  2018.
\newblock ISSN 1064-8275.

\bibitem[Estrin et~al.(2019{\natexlab{a}})Estrin, Friedlander, Orban, and
  Saunders]{EstrFrieOrbaSaun:2019a}
R.~Estrin, M.~P. Friedlander, D.~Orban, and M.~A. Saunders.
\newblock Implementing a smooth exact penalty function for equality-constrained
  nonlinear optimization.
\newblock \emph{{SIAM} J. Matrix Anal. Appl.}, (to appear), 2019{\natexlab{a}}.

\bibitem[Estrin et~al.(2019{\natexlab{b}})Estrin, Orban, and
  Saunders]{EstOrbSaun:2019LNLQ}
R.~Estrin, D.~Orban, and M.~A. Saunders.
\newblock {LNLQ}: An iterative method for linear least-norm problems with an
  error minimization property.
\newblock \emph{{SIAM} J. Matrix Anal. Appl.}, (to appear), 2019{\natexlab{b}}.

\bibitem[Fletcher(1970)]{Fletcher:1970}
R.~Fletcher.
\newblock A class of methods for nonlinear programming with termination and
  convergence properties.
\newblock In J.~Abadie, editor, \emph{Integer and Nonlinear Programming}, pages
  157--175. North-Holland, Amsterdam, 1970.

\bibitem[Fletcher(1973{\natexlab{a}})]{Fletcher:1973}
R.~Fletcher.
\newblock A class of methods for nonlinear programming: {III. Rates} of
  convergence.
\newblock In F.~A. Lootsma, editor, \emph{Numerical Methods for Nonlinear
  Optimization}. Academic Press, New York, 1973{\natexlab{a}}.

\bibitem[Fletcher(1973{\natexlab{b}})]{Fletcher:1973b}
R.~Fletcher.
\newblock An exact penalty function for nonlinear programming with
  inequalities.
\newblock \emph{Math. Programming}, 5:\penalty0 129--150, 1973{\natexlab{b}}.
\newblock ISSN 0025-5610.
\newblock \doi{10.1007/BF01580117}.

\bibitem[Freund and Nachtigal(1991)]{FreuNach:1991}
R.~W. Freund and N.~M. Nachtigal.
\newblock Q{MR}: a quasi-minimal residual method for non-{H}ermitian linear
  systems.
\newblock \emph{Numer. Math.}, 60\penalty0 (3):\penalty0 315--339, 1991.
\newblock ISSN 0029-599X.
\newblock \doi{10.1007/BF01385726}.
\newblock URL \url{https://doi.org/10.1007/BF01385726}.

\bibitem[Gersborg-Hansen et~al.(2006)Gersborg-Hansen, Bends\o{}e, and
  Sigmund]{GersBendSigm:2006}
A.~Gersborg-Hansen, M.~P. Bends\o{}e, and O.~Sigmund.
\newblock Topology optimization of heat conduction problems using the finite
  volume method.
\newblock \emph{Struct. Multidiscip. Optim.}, 31\penalty0 (4):\penalty0
  251--259, 2006.
\newblock ISSN 1615-147X.
\newblock \doi{10.1007/s00158-005-0584-3}.

\bibitem[Gould et~al.(2003)Gould, Orban, and Toint]{GoulOrbaToin:2003}
N.~I.~M. Gould, D.~Orban, and {\mbox{Ph}}.~L. Toint.
\newblock {CUTEr} and {SifDec}: A constrained and unconstrained testing
  environment, revisited.
\newblock \emph{{ACM} Trans. Math. Softw.}, 29\penalty0 (4):\penalty0 373--394,
  Dec. 2003.

\bibitem[Heinkenschloss and Ridzal(2014)]{HeinRidz:2014}
M.~Heinkenschloss and D.~Ridzal.
\newblock A matrix-free trust-region {SQP} method for equality constrained
  optimization.
\newblock \emph{SIAM J. Optim.}, 24\penalty0 (3):\penalty0 1507--1541, 2014.
\newblock \doi{10.1137/130921738}.

\bibitem[Leyffer(2006)]{Leyf:2006}
S.~Leyffer.
\newblock Complementarity constraints as nonlinear equations: theory and
  numerical experience.
\newblock In \emph{Optimization with Multivalued Mappings}, volume~2 of
  \emph{Springer Optim. Appl.}, pages 169--208. Springer, New York, 2006.
\newblock \doi{10.1007/0-387-34221-4_9}.

\bibitem[Maratos(1978)]{Mara:1978}
N.~Maratos.
\newblock \emph{Exact Penalty Function Algorithms for Finite Dimensional and
  Optimization Problems}.
\newblock PhD thesis, Imperial College of Science and Technology, London, UK,
  1978.

\bibitem[Nocedal and Wright(2006)]{NocedalW:2006}
J.~Nocedal and S.~J. Wright.
\newblock \emph{Numerical Optimization}.
\newblock Springer, New York, second edition, 2006.

\bibitem[Rees et~al.(2010)Rees, Dollar, and Wathen]{ReesDollWath:2010}
T.~Rees, H.~S. Dollar, and A.~J. Wathen.
\newblock Optimal solvers for pde-constrained optimization.
\newblock \emph{SIAM J. Sci. Comput.}, 32\penalty0 (1):\penalty0 271--298, Feb.
  2010.
\newblock ISSN 1064-8275.
\newblock \doi{10.1137/080727154}.

\bibitem[Ridzal(2013)]{Ridz:2013}
Ridzal.
\newblock Preconditioning of a full-space turst-region sqp algorithm for
  pde-constrained optimization.
\newblock \emph{Numerical Methods for PDE Constrained Optimization with
  Uncertain Data, Oberwolfach Reports}, 10\penalty0 (1):\penalty0 274--277,
  2013.

\bibitem[Saad and Schultz(1986)]{SaadSchu:1986}
Y.~Saad and M.~H. Schultz.
\newblock G{MRES}: a generalized minimal residual algorithm for solving
  nonsymmetric linear systems.
\newblock \emph{SIAM J. Sci. Statist. Comput.}, 7\penalty0 (3):\penalty0
  856--869, 1986.
\newblock \doi{10.1137/0907058}.

\bibitem[Sigmund and Petersson(1998)]{SigmPete:1998}
O.~Sigmund and J.~Petersson.
\newblock Numerical instabilities in topology optimization: A survey on
  procedures dealing with checkerboards, mesh-dependencies and local minima.
\newblock \emph{Structural optimization}, 16:\penalty0 68--75, 1998.

\bibitem[Steihaug(1983)]{Stei:1983}
T.~Steihaug.
\newblock The conjugate gradient method and trust regions in large scale
  optimization.
\newblock \emph{SIAM J. Numer. Anal.}, 20\penalty0 (3):\penalty0 626--637,
  1983.
\newblock \doi{10.1137/0720042}.

\bibitem[Stoll and Wathen(2012)]{StolWath:2012}
M.~Stoll and A.~Wathen.
\newblock Preconditioning for partial differential equation constrained
  optimization with control constraints.
\newblock \emph{Numer. Linear Algebra Appl.}, 19\penalty0 (1):\penalty0 53--71,
  2012.
\newblock ISSN 1070-5325.
\newblock \doi{10.1002/nla.823}.
\newblock URL \url{https://doi.org/10.1002/nla.823}.

\bibitem[Zavala and Anitescu(2014)]{AnitZava:2014}
V.~M. Zavala and M.~Anitescu.
\newblock Scalable nonlinear programming via exact differentiable penalty
  functions and trust-region {N}ewton methods.
\newblock \emph{SIAM J. Optim.}, 24\penalty0 (1):\penalty0 528--558, 2014.

\end{thebibliography}

% \newpage
% \tableofcontents
% \listoftodos\relax

\end{document}